\documentclass[reqno,centertags, 12pt]{amsart}
\usepackage{amsmath,amsthm,amscd,amssymb}
\usepackage{latexsym}
\newcommand{\bbC}{{\mathbb{C}}}
\newcommand{\bbD}{{\mathbb{D}}}

\newcommand{\bbZ}{{\mathbb{Z}}}

\newcommand{\calC}{{\mathcal{C}}}

\newcommand{\calH}{{\mathcal H}}

\newcommand{\calR}{{\mathcal R}}

\newcommand{\bdone}{{\boldsymbol{1}}}
\newcommand{\bddot}{{\boldsymbol{\cdot}}}

\newcommand{\lb}{\label}
\newcommand{\f}{\frac}

\newcommand{\ol}{\overline}
\newcommand{\ti}{\tilde  }

\newcommand{\dist}{\text{\rm{dist}}}

\newcommand{\spec}{\text{\rm{spec}}}

\newcommand{\ran}{\text{\rm{ran}}}

\newcommand{\ess}{\text{\rm{ess}}}

\newcommand{\s}{\text{\rm{s}}}

\newcommand{\supp}{\text{\rm{supp}}}

\newcommand{\bi}{\bibitem}

\newcommand{\beq}{\begin{equation}}
\newcommand{\eeq}{\end{equation}}
\newcommand{\ba}{\begin{align}}
\newcommand{\ea}{\end{align}}
\newcommand{\veps}{\varepsilon}

\newcounter{smalllist}
\newenvironment{SL}{\begin{list}{{\rm\roman{smalllist})}}{%
\setlength{\topsep}{0mm}\setlength{\parsep}{0mm}\setlength{\itemsep}{0mm}%
\setlength{\labelwidth}{2em}\setlength{\leftmargin}{2em}\usecounter{smalllist}%
}}{\end{list}}

\DeclareMathOperator{\Real}{Re}
\DeclareMathOperator{\Ima}{Im}

\allowdisplaybreaks
\numberwithin{equation}{section}

\newtheorem{theorem}{Theorem}[section]

\newtheorem*{p2.1}{Proposition 2.1}
\newtheorem{proposition}[theorem]{Proposition}
\newtheorem{lemma}[theorem]{Lemma}

\theoremstyle{definition}
\newtheorem{example}[theorem]{Example}

\theoremstyle{remark}
\newtheorem*{remark}{Remark}
\newtheorem*{remarks}{Remarks}

\newcommand{\abs}[1]{\lvert#1\rvert}

\title[Schur and Related Flows]{Zeros of OPUC and
Long Time Asymptotics of Schur and Related Flows}
\author[Barry Simon]{}


\email{bsimon@caltech.edu}
\thanks{$^*$ Supported in part by NSF grant DMS-0140592 and
U.S.--Israel Binational Science Foundation (BSF) Grant No.\ 2002068}

\begin{document}
\maketitle

\centerline{\scshape Barry Simon}
\medskip
{\footnotesize
 \centerline{Mathematics 253-37}
  \centerline{California Institute of Technology}
   \centerline{Pasadena, CA 91125, USA}
}

\medskip


\begin{abstract} We provide a complete analysis of the asymptotics for
the semi-infinite Schur flow: $\alpha_j(t)=(1-\abs{\alpha_j(t)}^2)
(\alpha_{j+1}(t)-\alpha_{j-1}(t))$ for $\alpha_{-1}(t)= 1$ boundary
conditions and $n=0,1,2,\dots$, with initial condition $\alpha_j(0)\in (-1,1)$.
We also provide examples with $\alpha_j(0)\in\bbD$ for which $\alpha_0(t)$
does not have a limit. The proofs depend on the solution via a direct/inverse
spectral transform.
\end{abstract}

\section{Introduction} \lb{s1}

One of our purposes in this paper is to study the long time asymptotics of the solution
of the differential equation on $(-1,1)^\infty$ for $t\geq 0$, $n\geq 0$,
\begin{equation} \lb{1.1}
\alpha'_n (t) = (1-\abs{\alpha_n(t)}^2) (\alpha_{n+1}(t) -\alpha_{n-1}(t))
\end{equation}
where $\alpha_{-1}$ is interpreted as
\begin{equation} \lb{1.2}
\alpha_{-1}(t) = -1
\end{equation}
for arbitrary boundary conditions $\alpha_n(0)\in (-1,1)$. This is called the Schur flow
\cite{AL75,AG94,FG99,Golppt}. We will also say something about complex initial conditions
with $\alpha_n(t)\in\bbD =\{z\mid\abs{z}<1\}$.

We got interested in this problem due to work of Golinskii \cite{Golppt} who proved for
initial conditions $\alpha_n(0)=0$ (all $n\geq 0$) that $\alpha_n(t)\to (-1)^n$,
and he obtains the leading $O(1/t)$ correction. From this point of view, our main
result is the following:

\begin{theorem}\lb{T1.1} Suppose each $\alpha_n(0)\in (-1,1)$. Then one of the following
holds:
\begin{SL}
\item[{\rm{(i)}}]
\begin{equation} \lb{1.3}
\alpha_n(t)\to (-1)^n \qquad\text{for all $n$}
\end{equation}
\item[{\rm{(ii)}}] There exists $1\leq N\leq \infty$ and
\begin{equation} \lb{1.4}
1>x_1 > x_2 > \cdots > x_N > -1
\end{equation}
so that
\begin{equation}\lb{1.5}
\begin{gathered}
\alpha_{2n}(t) \to 1 \quad n\geq 0, \qquad
\alpha_{2n-1}\to -x_n \quad \text{all } 1\leq n < N \\
\alpha_{2n-1} \to -x_N \quad \text{all } n\geq N
\end{gathered}
\end{equation}

\item[{\rm{(iii)}}] There exists $1\leq N \leq \infty$ and $x_j$ obeying \eqref{1.4} so that
\begin{equation}\lb{1.6}
\begin{gathered}
\alpha_{2n+1}(t) \to -1 \quad n\geq 0, \qquad
\alpha_{2n}(t)\to x_{n+1} \quad  0\leq n < N \\
\alpha_{2n}(t) \to x_N \quad  n\geq N -1
\end{gathered}
\end{equation}
\end{SL}
\end{theorem}

We will have a complete spectral theory analysis of which case one has based on
the initial conditions. But for now, we note:

\begin{proposition}\lb{P1.2}
\begin{SL}
\item[{\rm{(i)}}] If $\alpha_n(0)\to 0$ as $n\to\infty$, then we are in case
{\rm{(i)}} of Theorem~\ref{T1.1}.

\item[{\rm{(ii)}}] If $(-1)^n \alpha_n(0)\to \pm 1$ as $n\to\infty$, then we are in
case {\rm{(i)}} of Theorem~\ref{T1.1}.

\item[{\rm{(iii)}}] If $0<a<1$ and $\alpha_n(0)\to a$ as $n\to\infty$, then we are in
case {\rm{(ii)}} of Theorem~\ref{T1.1} and
\begin{equation} \lb{1.7}
x_N = 1-2a^2 \text{ if } N<\infty \quad{\text{or}}\quad x_n\downarrow 1-2a^2
\text{ if } N=\infty
\end{equation}

\item[{\rm{(iv)}}] If $-1 <a<0$ and $\alpha_n(0)\to a$ as $n\to\infty$, then we are in
case {\rm{(iii)}} of Theorem~\ref{T1.1} and \eqref{1.7} holds.

\item[{\rm{(v)}}] Case {\rm{(ii)}} holds with $x_n\downarrow -1$ as $n\to\infty$ if and only
if $\alpha_n(0)\to 1$ as $n\to\infty$.

\item[{\rm{(vi)}}] Case {\rm{(iii)}} holds with $x_n\downarrow -1$ as $n\to\infty$ if and only
if $\alpha_n(0)\to -1$ as $n\to\infty$.
\end{SL}
\end{proposition}

{\it Note.} These possibilities are consistent with the partial results of Theorem~5
of Golinskii \cite{Golppt}.

\smallskip
The situation is rather more subtle if we allow complex initial conditions:

\begin{proposition}\lb{P1.3} There exists $\{\alpha_n(0)\}\in\bbD^\infty$ so that
$\alpha_0 (t)$ does not have a limit.
\end{proposition}

Besides this, we will discuss rates of convergence. In cases (ii) and (iii) for
$j<2N+1$ (resp.\ $j<2N$), the rate will be exponentially fast. In other cases,
the situation can be subtle, although if $\sum_{n=0}^\infty \abs{\alpha_n(0)} <
\infty$, the rate is that found by Golinskii when $\alpha_n(0)=0$.

At first glance, it seems surprising that one can obtain such detailed information
for a nonlinear equation. The reason, of course, is that \eqref{1.1} is completely
integrable. Indeed, it is exactly solvable via a spectral transform \cite{FG99,KN2,Golppt}.
In this sense, this problem is a close analog of work of Moser \cite{Mo1,Mo2} and
Deift--Li--Tomei \cite{DLT85} on asymptotics of Toda flows, except for an extra subtlety
we will discuss shortly.

Just as Toda is closely connected to the theory of orthogonal polynomials on
the real line (OPRL), the theory of orthogonal polynomials on the unit circle
(OPUC) \cite{Szb,GBk1,OPUC1,OPUC2,1Foot} will be central here. As Golinskii \cite{Golppt}
notes for Toda, OPRL ``plays one of the first fiddles in the performance
(albeit not entering the final results directly)." To push his metaphor, the present
paper promotes OPUC and OPRL to concert soloist---OPs enter directly into the results
(see Theorem~\ref{T1.5} below), and more directly in our proofs than in previous works.

Recall (e.g., \cite{1Foot,OPUC1}) that nontrivial probability measures,
$d\mu$, on $\partial\bbD =\{z\mid \abs{z}=1\}$ are parametrized by
$\{\alpha_n\}_{n=0}^\infty\in\bbD^\infty$ via the Szeg\H{o} recursion relations
\begin{align}
\Phi_{n+1}(z) &= z\Phi_n(z) -\bar\alpha_n \Phi_n^*(z) \lb{1.8} \\
\Phi_n^*(z) &= z^n\, \ol{\Phi_n (1/\bar z)}\lb{1.9}
\end{align}
where $\Phi_n(z)$ are the monic orthogonal polynomials in $L^2 (\partial\bbD,d\mu)$.
We often use $\alpha_n (d\mu)$ when we want to make the $d\mu$-dependence explicit.
The situation of \eqref{1.1} is made explicit by

\begin{proposition}[\cite{FG99,Golppt}; also see our first appendix and \cite{CSprep}]\lb{P1.4}
Given $\alpha_n(0)\in\bbD^\infty$, define $d\mu$ by
\begin{equation} \lb{1.10}
\alpha_n(d\mu)= \alpha_n(0)
\end{equation}
and $d\mu_t$ by
\begin{equation} \lb{1.11}
d\mu_t(\theta)= \f{e^{2t\cos(\theta)}\, d\mu(\theta)}{\int e^{2t\cos (\theta)}\, d\mu(\theta)}
\end{equation}
Then
\begin{equation} \lb{1.12}
\alpha_n(t) \equiv \alpha_n (d\mu_t)
\end{equation}
is the unique solution of \eqref{1.1} with initial conditions $\alpha_n(0)$.
\end{proposition}

Following the analog of Deift--Li--Tomei for Toda flows \cite{DLT85}, we will want to consider
generalized Schur flows associated to any bounded real-valued function $G(\theta)$
on $\partial\bbD$ via
\begin{align}
\Sigma_t (d\mu)(\theta) &= \f{e^{tG(\theta)}\, d\mu(\theta)}
{\int e^{tG(\theta)}\, d\mu(\theta)} \lb{1.13} \\
\sigma_t (\alpha_n(0)) &= \alpha_n (\Sigma_t (d\mu)) \lb{1.14}
\end{align}
When we want to make $G$ explicit, we will use $\sigma_t^G$ and $\Sigma_t^G$. We will
say more about these flows in the first appendix. \cite{KN2} calls these generalized
Ablowitz--Ladik flows. Since generalized Schur flows preserve reality only if
$G(-\theta)=G(\theta)$, there is some reason in their choice, but we prefer to
emphasize the connection to Schur functions and OPUC.

One can also define $\Sigma_t$ and $\sigma_t$ for trivial measures, that is,
measures on $\partial\bbD$ with finite support
\begin{equation} \lb{1.15}
d\mu = \sum_{j=1}^m \mu_j \delta_{e^{i\theta_j}}
\end{equation}
parametrized by $\{\alpha_j\}_{j=0}^{m-1}$ with $\alpha_0, \dots, \alpha_{m-2}\in
\bbD$ and $\alpha_{m-1} =(-1)^{m+1} \prod_{j=1}^m e^{i\theta_j}$. In this case,
since all that matters are $\{G(\theta_j)\}_{j=1}^m$, we can suppose $G$ is a
polynomial. For this case, the asymptotics of generalized Schur flows were studied
by Killip--Nenciu \cite{KN2}.

We will also analyze the long time asymptotics in this case. To do so, we will let
\begin{equation} \lb{1.16}
z_j = e^{i\theta_j}
\end{equation}
and renumber, so
\begin{equation} \lb{1.17}
G(z_1)\geq G(z_2)\geq \cdots \geq G(z_m)
\end{equation}
We define the $K$-groups $K_1, \dots, K_\ell$ to be those indices $K_1\equiv \{1,
\dots, k_1\}$, $K_2\equiv \{k_1 +1, \dots k_2\}, \dots$, $K_\ell =
\{k_{\ell-1}+1, \dots, k_\ell\equiv m\}$, so $G(z_j) =G(z_p)$ if $j,p\in K_\ell$
and so $G(z_{k_j}) > G(z_{k_j+1})$. Thus, the $K$-groups are the level sets of
$G$ on $\{z_j\}_{j=1}^m$.

It will also be convenient to define
\begin{equation} \lb{1.18}
z^{(K_j)} = \prod_{p\leq k_{j-1}} z_j
\end{equation}
the product of those $z$ where $G$ is larger than the value common in $K_j$.
Also, given an initial point measure, $d\mu$ of the form \eqref{1.15},
we define the $K$-group induced measure by
\begin{equation} \lb{1.19}
d\mu^{(K_j)} =\sum_{\ell=k_{j-1}+1}^{k_j} \mu_\ell^{(K_j)} \delta_{e^{i\theta_\ell}}
\end{equation}
where
\begin{equation} \lb{1.20}
\mu_\ell^{(K_j)} = \f{[\prod_{p=1}^{k_{j-1}} \abs{z_\ell - z_p}^2 ]\mu_\ell}
{\sum_{m=k_{j-1}+1}^{k_j}\, [\prod_{p=1}^{k_{j-1}} \abs{z_m - z_p}^2]\mu_m}
\end{equation}
Then we will prove:

\begin{theorem}\lb{T1.5} Let $d\mu$ be a finite measure given by \eqref{1.15} and $G$
defined on $\{z_j\}_{j=1}^m$ and real-valued. Then for $\ell\in K_j$,
\begin{equation} \lb{1.21}
\sigma_t (\alpha_{\ell-1}) \to (-1)^{k_{j-1}}\, \ol{z^{(K_j)}}\,
\alpha_{\ell-k_{j-1}-1} (d\mu^{(K_j)})
\end{equation}
\end{theorem}

Killip--Nenciu also obtain a limit theorem; we will discuss its relation to ours in
Section~\ref{s3}.

At first sight, the problem looks very easy. $\alpha_n (\cdot)$ is continuous under
weak convergence (a.k.a.\ vague or weak-$*$) convergence of measures, so one need only
find the weak limit of $d\mu_t$ or $\Sigma_t (d\mu)$. For the Schur flow case, $d\mu_t$
has either a one- or two-point support for its weak limit. But for trivial (i.e.,
finite point) limits, $\lim \alpha_n(\cdot)$ is only determined for $n$ smaller
than the number of points in the support of the limit, so this only determines at
most two $\alpha$'s! The issue can be seen clearly in the context of Theorem~\ref{T1.5}:
$d\mu_t$ converges to $d\mu^{(K_1)}$ and so only determines $\{\alpha_{\ell-1}
(\infty)\}_{\ell=1}^{k_1}$.

The key to the general analysis is to track the zeros, $\{z_j^{(n)}\}_{j=1}^n$,
of the OPUC, $\Phi_n(z)$. By \eqref{1.8}, $\Phi_n^*(0) =1$ and
\begin{equation} \lb{1.22}
\Phi_n(z) =\prod_{j=1}^n (z-z_j^{(n)})
\end{equation}
so we have
\begin{equation} \lb{1.23}
\alpha_n = (-1)^n \prod_{j=1}^{n+1}\, \ol{z_j^{(n+1)}}
\end{equation}
The zeros in turn are determined by the Szeg\H{o} variational principle, that
\begin{equation} \lb{1.24}
\int \biggl|\, \prod_{j=1}^n (z-w_j)\biggr|^2\, d\mu(\theta)
\end{equation}
is minimized precisely with $\{w_j\}_{j=1}^n =\{z_j^{(n)}\}_{j=1}^n$.

We will find a consequence of the variational principle (Theorem~\ref{T3.3A}
below) that lets us use weak convergence beyond the naive limit.

These same ideas work for the Toda flows and are in some ways simpler
there since the exponential factor there, $e^{2tx}$, is strictly monotone.
We begin in Section~\ref{s2}, as a warmup, by proving Moser's theorem on the
asymptotics for finite Toda flows. In Section~\ref{s3}, we prove Theorem~\ref{T1.5}.
In Section~\ref{s4}, we prove Theorem~\ref{T1.1} and Proposition~\ref{P1.2}---as in
Deift--Li--Tomei, the extreme points in the essential spectrum are crucial.
Section~\ref{s5} discusses second-order corrections. In Section~\ref{s6}, we
present the example of Proposition~\ref{P1.3}, and in Section~\ref{s7},
discuss hypotheses which prevent the pathologies there. Appendix~A combines
ideas of Deift--Li--Tomei \cite{DLT85,DNT83} for Toda flows with notions from
Killip--Nenciu \cite{KN2} to talk about difference equations associated
with \eqref{1.13}. Appendix~B provides a new proof and strengthening of
results of Denisov--Simon on zeros of OPUC near isolated points of $\supp (d\mu)$.

We should note that while we show the asymptotics of $\alpha_n(t)$ for $n$ fixed
and $t\to\infty$ is simple, and one can see that, often, asymptotics of $\alpha_n(t)$
for $t$ fixed and $n\to\infty$ is easy (e.g., Golinskii \cite{Golppt} shows that
$\alpha_n(0)\in\ell^1$ (resp.\ $\ell^2$) implies $\alpha_n(t)\in \ell^1$ (resp.\ $\ell^2$)),
the subtle asymptotics is for $\alpha_n(t)$ as $t\to\infty$ and $n/t\to q\in (0,\infty)$.
For $\alpha_n(0)\equiv 0$, this is studied using Riemann--Hilbert methods in \cite{BDJ99},
and no doubt their methods extend to any case $\abs{\alpha_n(0)}\leq e^{-cn}$ for $c>0$.
Indeed, using ideas of \cite{MMppt}, one can probably handle other classes where $d\mu$ is
not analytic.

\section{Moser's Theorem on Toda Asymptotics} \lb{s2}

As a warmup, we consider a probability measure
\begin{equation} \lb{2.1}
d\rho(x) =\sum_{j=1}^N \rho_j \delta_{x_j}
\end{equation}
where each $\rho_j >0$, and the family of measures
\begin{equation} \lb{2.2}
d\rho_t (x) = \f{e^{2tx}\,d\rho(x)}{\int e^{2tx}\,d\rho(x)}
\end{equation}
The Jacobi parameters, $\{a_n\}_{n=1}^{N-1}\cup\{b_n\}_{n=1}^N$, associated to $d\rho$
are defined by looking at the recursion relations associated to the monic orthogonal
polynomials
\begin{equation} \lb{2.2a}
xP_n(x) =P_{n+1}(x) + b_{n+1} P_n(x) + a_n^2 P_{n-1}(x)
\end{equation}
with $a_n >0$. The Jacobi parameters associated to $d\rho_t$ obey the Toda equations in
Flaschka form:
\begin{align}
\f{da_n}{dt} &=a_n (b_{n+1}-b_n)  \lb{2.3}\\
\f{db_n}{dt} &= 2(a_n^2 - a_{n-1}^2) \lb{2.4}
\end{align}
(with $a_0=a_N\equiv 0$ in \eqref{2.4}).

We order the $x_j$'s by
\[
x_1 > x_2 > \cdots > x_N
\]
Then our main result in this section is the following:

\begin{theorem}\lb{T2.1} The above finite Toda chains obey
\begin{equation} \lb{2.5}
\lim_{t\to\infty}\, b_j(t) = x_j \qquad
\lim_{t\to\infty}\, a_j(t) =0
\end{equation}
{\rm{(}}the first for $j=1,\dots,N$ and the second for $j=1,2,\dots, N-1${\rm{)}}.
Moreover, errors for $b_n$ are $O(e^{-ct})$ and for $a_n$ are $O(e^{-ct/2})$ where
\begin{equation} \lb{2.6}
c=\min_{j=1,2,\dots, N-1} (x_j-x_{j+1})
\end{equation}
\end{theorem}

\begin{remarks} 1. This is a celebrated result of Moser \cite{Mo1}, but it is interesting
to see it proven using zeros of OPRL. In any event, it is a suitable warmup for our result
on OPUC.

\smallskip
2. The $O(e^{-ct})$ and $O(e^{-ct/2})$ estimates are not ideal; we will discuss this further
and obtain finer error estimates at the end of this section.

\smallskip
3. The same proof shows that as $t\to -\infty$, $b_j\to x_{N+1-j}$ and $a_j\to 0$.
\end{remarks}

$P_j(x)$ has $j$ simple zeros $x_1^{(j)} >x_2^{(j)} >\cdots > x_j^{(j)}$. (It is known
that $x_k > x_k^{(j+1)} > x_k^{(j)} > x_{n-j+k-1}$, but we won't need that.) The key fact
is:

\begin{theorem}\lb{T2.2} For each $j=1, \dots, N$ and $k=1, \dots, j$,
\begin{equation} \lb{2.6x}
\lim_{t\to\infty}\, x_k^{(j)}(t) = x_k
\end{equation}
The errors are $O(e^{-ct})$ with $c$ given by \eqref{2.6}.
\end{theorem}

\begin{remark} Since $P_N(x)=\prod_{j=1}^N (x-x_j)$, we have $x_j^{(N)}=x_j$ for all $t$.
\end{remark}

\begin{proof} $P_j(x)$ is the projection of $x^j$ on $\{1, \dots, x^{j-1}\}^\perp$, so for any
monic polynomial $Q(x)$ of degree $j$,
\begin{equation} \lb{2.7}
\int \abs{Q(x)}^2\, d\rho(x) \geq \int \abs{P_j(x)}^2\, d\rho(x)
\end{equation}
Pick $Q(x)=\prod_{\ell=1}^j (x-x_\ell)$, which minimizes the contributions of $x_1, \dots, x_j$
to the integral, and see that
\begin{equation} \lb{2.8}
\int \abs{P_j(x)}^2\, d\rho(x) \leq e^{2tx_{j+1}} (x_1 -x_N)^{2j}
\end{equation}
since $\abs{Q(x_\ell)}\leq \abs{x_\ell -x_1}^j$ for $\ell\geq j+1$. On the other hand,
\begin{equation} \lb{2.8a}
\abs{P_j(x)}\geq \bigl[\, \min_{\ell=1, \dots, j} \, \abs{x-x_\ell^{(j)}}\bigr]^j
\end{equation}
so for $q=1, \dots, j$,
\begin{equation} \lb{2.9}
\rho_q e^{2tx_q} \min_{\ell=1, \dots, j}\, [\abs{x_q - x_\ell^{(j)}}^j ]^2 \leq
\int \abs{P_j(x)}^2\, d\rho(x)
\end{equation}
We conclude for $q\leq j$ that
\begin{equation} \lb{2.10}
\min_{\ell=1, \dots, j}\, \abs{x_q -x_\ell^{(j)}}^{2j} \leq
\rho_q^{-1} (x_1-x_N)^{2j} e^{2t (x_{j+1}-x_q)}
\end{equation}
which shows that each $x_q$ has an $x_\ell^{(j)}$ exponentially near to it, but only
$O(e^{-ct/j})$.

But once we know each such $x_q$ has one zero exponentially near, we see that for $t$
large, all other zeros are a distance at least $\f12 c$ away. Thus \eqref{2.8a} can be
replaced, for $t$ large, by
\begin{equation} \lb{2.11}
\abs{P_j(x)} \geq \bigg(\f{c}{2}\biggr)^{j-1} \min_{\ell=1, \dots, j}\, \abs{x-x_\ell^{(j)}}
\end{equation}
Plugging this into \eqref{2.9} and finding the analog of \eqref{2.10} leads to an
$O(e^{-ct})$ error. Explicitly, \eqref{2.10} is replaced by
\begin{equation}\lb{2.14a}
\min_{\ell=1, \dots, j}\, \abs{x_q-x_\ell^{(j)}}^2 \leq C e^{2t(x_{j+1}-x_q)}
\end{equation}
\end{proof}

\begin{remark} We will need a better error estimate in the next section
and show how to get it later in the section.
\end{remark}

\begin{proof}[Proof of Theorem~\ref{T2.1}] \eqref{2.2a} can be rewritten:
\begin{equation} \lb{2.12}
\prod_{\ell=q}^{n+1} (x-x_j^{(n+1)}) = (x-b_{n+1}) \prod_{\ell=1}^n (x-x_j^{(n)})
- a_n^2 P_{n-1}(x)
\end{equation}
Identifying the $x^n$ and $x^{n-1}$ terms, we see the analog of \eqref{1.23};
the first for $n=1,2,\dots, N$ and the second for $n=1,2,\dots, N-1$:
\begin{align}
b_{n+1} &= \sum_{j=1}^{n+1} x_j^{(n+1)} - \sum_{j=1}^n x_j^{(n)} \lb{2.13} \\
a_n^2 &= b_{n+1} \sum_{j=1}^n x_j^{(n)} + \sum_{1\leq j < \ell \leq n} x_j^{(n)}
x_\ell^{(n)} - \sum_{1\leq j <\ell \leq n+1} x_j^{(n+1)} x_\ell^{(n+1)} \lb{2.14}
\end{align}
\eqref{2.6x} and the error estimates of Theorem~\ref{T2.2} immediately imply
$b_j(t)-x_j = O(e^{-ct})$ and $a_j(t)^2 = O(e^{-ct})$.
\end{proof}

We next want to note, following Moser, that the differential equations
\eqref{2.3} and \eqref{2.4} yield better error estimates than Theorem~\ref{T2.1}
has and then explain how to improve the estimates on zeros to get better estimates
on the errors of $b_n$ and $a_n$ with the zeros.

Once we know $b_j(t)\to x_j$, \eqref{2.2} implies
\begin{equation} \lb{2.15}
t^{-1} \log a_j(t) \to x_{j+1} - x_j <0
\end{equation}
Indeed, since the approach of $b_j(t)$ to $x_j$ is exponentially fast,
\begin{equation} \lb{2.16}
\log a_j(t)  -t(x_{j+1} -x_j) \to \log C_j
\end{equation}
for some finite $C_j$, and thus,
\begin{equation} \lb{2.17}
a_j\sim C_j e^{-t(x_{j+1}-x_j)}
\end{equation}
proving the error should be $O(e^{-ct})$, not $O(e^{-ct/2})$. Then plugging \eqref{2.17}
into \eqref{2.4}, we see that
\begin{equation} \lb{2.18}
\abs{b_j(t)-x_j} \leq \ti C_j \exp (-2t \min[(x_{j+1} -x_j), (x_j-x_{j-1})])
\end{equation}
with the right side being the exact order of error if $x_{j+1}-x_j \neq x_j - x_{j-1}$
(if there is equality, $a_j^2$ and $a_{j-1}^2$ can completely or partially cancel).
Thus, the error is $O(e^{-2ct})$, not $O(e^{-ct})$.

To improve our estimates on zeros, we use the minimization principle to get a
self-consistency equation on the zeros. This result, proven using orthogonality of
$P$ to $P/(x-x_0)$ is well-known (see (3.3.3) of \cite{Szb}); we give a variational
principle argument in line with the strategy in this paper. The OPUC analog
is (1.7.51) of \cite{OPUC1}.

\begin{lemma}\lb{L2.4} The zeros $x_k^{(j)}$ of $P_j(x)$ obey
\begin{equation} \lb{2.20}
x_k^{(j)}= \f{\int x\prod_{\ell\neq k} \abs{x-x_\ell^{(j)}}^2\, d\rho(x)}
{\int \prod_{\ell\neq k} \abs{x-x_\ell^{(j)}}^2\, d\rho(x)}
\end{equation}
In particular, for any $y$,
\begin{equation} \lb{2.21}
\abs{x_k^{(j)}-y}\leq \f{\int \abs{x-y} \prod_{\ell\neq k} \abs{x-x_\ell^{(j)}}^2\, d\rho(x)}
{\int \prod_{\ell\neq k} \abs{x-x_\ell^{(j)}}^2\, d\rho(x)}
\end{equation}
\end{lemma}

\begin{proof} Since $\int \prod_{n=1}^j \abs{x-y_\ell}^2\, d\rho(x)$ is minimized at $y_\ell
=x_\ell^{(j)}$, the derivative with respect to $y_k$ at this point is zero, that is,
\[
\int (x-x_k^{(j)}) \prod_{\ell\neq k} \abs{x-x_\ell^{(j)}}\, d\rho(x)=0
\]
which is \eqref{2.20}. \eqref{2.20} implies \eqref{2.21} by noting $x_k^{(j)}-y$ is given
by \eqref{2.20} with the first $x$ in the integrand replaced by $x-y$.
\end{proof}

\begin{theorem} \lb{T2.5} For $j=1, \dots, N$ and $k=1, \dots, j$,
\begin{equation} \lb{2.22}
\abs{x_k^{(j)}(t)-x_k} \leq Ce^{-2t(x_k-x_{j+1})}
\end{equation}
\end{theorem}

\begin{proof} We begin by noting that since the $\{x_\ell^{(j)}(t)\}_{\ell\neq k}$ for
$t$ large are very near $x_\ell$ (and so, not near $x_k$) that for some $T_0$ and constant
$C_1$ and all $t\geq T_0$,
\begin{equation} \lb{2.23}
\int \prod_{\ell\neq k}\, \abs{x-x_\ell^{(j)}}^2\, d\rho_t(x) \geq C_1 e^{2tx_k}
\end{equation}
Moreover, since $x_1 > x_2 > \dots$,
\begin{equation} \lb{2.24}
\sum_{m\geq j+1} \rho_m e^{2tx_m} \abs{x_m-x_k} \prod_{\ell\neq k} \, \abs{x_m-x_\ell^{(j)}}^2
\leq Ce^{2t x_{j+1}}
\end{equation}
so we need only control the terms $m=1,2,\dots, j$ in estimating \eqref{2.21}.

For $m=1, \dots, j$, we use \eqref{2.14a} to see
\begin{align*}
\rho_n e^{2tx_m} \abs{x_m-x_k}\, \prod_{\ell\neq k}\, \abs{x_m-x_\ell^{(j)}}^2
&\leq C e^{2tx_m} e^{2t(x_{j+1} -x_m)} \\
&\leq C e^{2tx_{j+1}}
\end{align*}
Thus all terms in the numerator of \eqref{2.21} with $y=x_k$ are bounded by $Ce^{2tx_{j+1}}$.
Combining this with \eqref{2.23}, we obtain \eqref{2.22}.
\end{proof}

\begin{remark} Putting the improved bound \eqref{2.22} in place of \eqref{2.14a} shows
that the sum in \eqref{2.24} dominates the sum in the numerator of \eqref{2.21}.
\end{remark}

\begin{theorem}\lb{T2.6} We have
\begin{align}
\abs{b_j(t)-x_j} &\leq C_3 [\exp (2t (x_{j+1}-x_j)) + \exp (2t (x_j-x_{j-1}))] \lb{2.27}\\
\abs{a_j(t)} &\leq C_4 [\exp (t(x_{j+2}-x_{j+1})) + \exp (t (x_{j+1}-x_j)) \notag \\
&\qquad \qquad \qquad + \exp(t(x_j-x_{j-1}))] \lb{2.28}
\end{align}
\end{theorem}

\begin{remark} As explained above (see \eqref{2.17} and \eqref{2.18}), \eqref{2.27}
is optimal, while \eqref{2.28} is not quite, although it has the proper $e^{-tc}$
behavior.
\end{remark}

\begin{proof} \eqref{2.27} follows from \eqref{2.13} and \eqref{2.22}, while \eqref{2.28}
follows from \eqref{2.14} and \eqref{2.22}.
\end{proof}

\section{A Theorem of Killip and Nenciu} \lb{s3}

In this section, we want to prove Theorem~\ref{T1.5}. We will follow the strategy of the
last section with some changes necessitated by the fact the $K$-groups can have more
than one point. In particular, we cannot use mere counting to be sure only one zero
approaches a single pure point. Instead we will need the following theorem of
Denisov--Simon that appears as Theorem~1.7.20 of \cite{OPUC1}:

\begin{theorem}\lb{T3.1} Let $z_0$ be an isolated point of the support of a probability
measure on $\partial\bbD$. Let
\begin{equation} \lb{3.1}
d=\dist (z_0, \supp(d\mu)\backslash\{z_0\})
\end{equation}
Then each OPUC, $\Phi_j(z;d\mu)$ has at most one zero in the circle of radius $d^2/6$
about $z_0$.
\end{theorem}

\begin{remarks} 1. For OPRL, isolated points of the support can have two nearby zeros;
see \cite{DS2}.

\smallskip
2. See Appendix~B for an alternate proof (of a stronger result, namely, getting
$d^2/4$ rather than $d^2/6$) that uses operator theory.
\end{remarks}

\begin{lemma}\lb{L3.2} The zeros $z_k^{(n)}$ of $\Phi_n (z;d\mu)$ for any $d\mu$ on
the unit circle obey
\begin{equation} \lb{3.2}
z_k^{(n)}= \f{\int z \prod_{\ell\neq k} \abs{z-z_\ell^{(n)}}^2\, d\mu(\theta)}
{\int \prod_{\ell\neq k} \abs{z-z_\ell^{(n)}}^2\, d\mu(\theta)}
\end{equation}
In particular, for any $y$,
\begin{equation} \lb{3.3}
\abs{z_k^{(n)}-y} \leq \f{\int \abs{z-y} \prod_{\ell\neq k} \abs{z-z_\ell^{(n)}}^2\, d\mu(\theta)}
{\int \prod_{\ell\neq k} \abs{z-z_\ell^{(n)}}^2\, d\mu(\theta)}
\end{equation}
\end{lemma}

\begin{remarks} 1. If $d\mu$ is trivial with $N$ points in its support, we need $n\leq N$.

\smallskip
2. In the integrals, $z=e^{i\theta}$.

\smallskip
3. \eqref{3.2} is (1.7.51) of \cite{OPUC1}. Again we give a variational proof.
\end{remarks}

\begin{proof} In $\int \prod_{\ell=1}^n \abs{z-y_\ell}^2\, d\mu(\theta)$, which is minimized
by $y_\ell=z_\ell^{(n)}$, the $y_\ell$'s are complex so we can write it as a function of
$y_\ell$ and $\bar y_\ell$ and demand all $\partial/\partial\bar y_\ell$ and $\partial/
\partial y_\ell$ vanish at $y_\ell = z_\ell^{(n)}$. \eqref{3.2} comes from the $\partial/
\partial\bar y_k$ derivative (or conjugate of the $\partial/\partial y_k$ derivative).
\eqref{3.3} follows immediately from \eqref{3.2}.
\end{proof}

\begin{theorem}\lb{T3.3} Let $\{z_j^{(n)}(t)\}_{j=1}^n$ be the zeros of $\Phi_n
(z;d\mu_t)$ where $n\in K_m$, $\mu_t$ is given by \eqref{1.13}, and $d\mu_{t=0}$ has
finite support. Here $K_m$ are the $K$-group defined after \eqref{1.17}. Then for
$t$ large, $\Phi_n$ has exactly one zero near each $\{z_j\}_{j=1}^{k_{m-1}}$ which,
by renumbering, we can suppose are $z_j^{(n)}(t)$. Moreover, for $1\leq j\leq k_{m-1}$,
\begin{equation} \lb{3.4}
\abs{z_j^{(n)}(t)-z_j} \leq C \exp (t [G(z_n)-G(z_j)])
\end{equation}
\end{theorem}

\begin{proof} The proof is identical to Theorem~\ref{T2.5}, given Theorem~\ref{T3.1}
to be sure $\prod_{\ell\neq j}\abs{z-z_\ell^{(n)}}^2$ stays away from the
zero at $z=z_j$.
\end{proof}

To continue, we will need a lovely consequence of the Szeg\H{o} variational
principle that will also be the key to the arguments in Section~\ref{s4}. Recall
that one can define monic OPUC for any positive measure, even if not normalized,
and, of course,
\begin{equation}\lb{3.4a}
\Phi_j (z;c\,d\mu) =\Phi_j (z;d\mu)
\end{equation}
for any positive constant $c$.

\begin{theorem}\lb{T3.3A} Let $d\mu$ be a nontrivial measure on $\partial\bbD$ and let
$\{z_j\}_{j=1}^k$ be among the zeros of $\Phi_n (z;d\mu)$. Then
\begin{equation}\lb{3.4b}
\Phi_n (z;d\mu) = \prod_{j=1}^k (z-z_j) \Phi_{n-k} \biggl( z; \prod_{j=1}^k \,
\abs{z-z_j}^2\, d\mu\biggr)
\end{equation}
\end{theorem}

\begin{remark} The $z_j$'s can be repeated up to their multiplicity.
\end{remark}

\begin{proof} Let $Q_{n-k}$ be a monic polynomial of degree $n-k$, so $\prod_{j=1}^k
(z-z_j) Q_{n-k}$ is a monic polynomial of degree $n$. Thus, by the Szeg\H{o}
variational principle,
\begin{align*}
\int \biggl| \f{\Phi_n (z;d\mu)}{\prod_{j=1}^k (z-z_j)}\biggr|^2 \prod_{j=1}^k \,
\abs{z-z_j}^2\, d\mu &= \int \abs{\Phi_n (z;d\mu)}^2\, d\mu \\
&\leq \int \biggl|\, \prod_{j=1}^k (z-z_j)Q_{n-k}(z)\biggr|^2\, d\mu \\
&= \int \abs{Q_{n-k}(z)}^2 \prod_{j=1}^k \, \abs{z-z_j}^2\, d\mu
\end{align*}

Since $\Phi_n (z;d\mu)/\prod_{j=1}^k (z-z_j)$ is a monic polynomial of degree
$n-k$ and $Q$ is arbitrary, the Szeg\H{o} variational principle implies that
\begin{equation}\lb{3.4c}
\f{\Phi_n(z;d\mu)}{\prod_{j=1}^k (z-z_j)} = \Phi_{n-k} \biggl( z;
\prod_{j=1}^k \, \abs{z-z_j}^2\, d\mu \biggr)
\end{equation}
which is \eqref{3.4b}.
\end{proof}

\begin{remarks} 1. One can also prove this using orthogonality. For if $\ell < n-k$,
$z^\ell \prod_{j=1}^k (z-z_j) \perp \Phi_n (z;d\mu)$, so
\begin{equation}\lb{3.4d}
\int \f{\Phi_n (z;d\mu)}{\prod_{j=1}^k \abs{z-z_j}^2}\,
\ol{z^\ell \prod_{j=1}^k (z-z_j)}\, \prod_{j=1}^k\, \abs{z-z_j}^2\, d\mu =0
\end{equation}
but
\[
\text{LHS of \eqref{3.4d}} =\int \ol{z^\ell}\, \f{\Phi_n (z;d\mu)}{\prod_{j=1}^k (z-z_j)}\,
\prod_{j=1}^k\, \abs{z-z_j}^2\, d\mu
\]
proving \eqref{3.4c}.

\smallskip
2. \eqref{3.4c} for $n-k=1$ is easily seen to be equivalent to \eqref{3.2}.
\end{remarks}

\begin{theorem}\lb{T3.4} Under the hypotheses of Theorem~\ref{T1.5} with
$d\mu_t =\Sigma_t (d\mu)$ and $\ell\in K_j$,
\begin{equation} \lb{3.5}
\lim_{t\to\infty}\, \Phi_\ell (z;d\mu_t) = \prod_{p=1}^{k_{j-1}} (z-z_p)
\Phi_{\ell-k_{j-1}-1} (z;d\mu^{(K_j)})
\end{equation}
\end{theorem}

\begin{proof} Let $\{z_p(t)\}_{p=1}^{k_{j-1}}$ be the zeros of $\Phi_\ell (z;d\mu_t)$
which converge to $\{z_p\}_{p=1}^{k_{j-1}}$ as $t\to\infty$. By \eqref{3.4b} and
\eqref{3.4a} for $n\in K_j$,
\begin{equation}\lb{3.6}
\Phi_\ell (z;d\mu_t) = \prod_{p=1}^{k_{j-1}} (z-z_p(t)) \Phi_{\ell-k_{j-1}}
\biggl( z;e^{-tG(z_n)} \prod_{p=1}^{k_{j-1}} \, \abs{z-z_p(t)}^2\, d\mu_t(z)\biggr)
\end{equation}
On account of \eqref{3.4}, the weights of $\{z_p\}_{p=1}^{k_{j-1}}$ in the measure on
the right in \eqref{3.6} are bounded by
\[
e^{2t(G(z_n)-G(z_y))} e^{-tG(z_n)} e^{+tG(z_j)}\to 0
\]
since $G(z_j)< G(z_n)$. The weights of points in $K_{j+1}, \dots, K_\ell$ go to zero.
Thus
\[
e^{-tG(z_n)} \prod_{p=1}^{k_{j-1}}\, \abs{z-z_p(t)}^2\, d\mu_t \to C\, d\mu^{(K_j)}
\]
where $C$ is a constant and the convergence is weak. For we have shown the contributions
of $z_\ell \notin K_j$ go to zero and the weights at $z_\ell\in K_j$ converge to
$\prod_{p=1}^{k_{j-1}} \abs{z_\ell - z_p}^2 \mu(\{z_\ell\})$ since $G(z_n)=G(z_\ell)$
and $z_p (t)\to z_p$.

\eqref{3.5} is immediate by continuity of OPs for index less or equal to the
number of points in the support of the limiting measure.
\end{proof}

\begin{remark} As stated, this theorem required $d\mu_{t=0}$ (and so $d\mu_t$) have finite
support. However, the proof works without change if $d\mu_{t=0}=d\nu_1 + d\nu_2$ with
$d\nu_1$ finite and $d\nu_2$ such that
\[
\sup_{z\in\supp(d\nu_2)} G(z) < \min_{z\in\supp(d\nu_1)} G(z)
\]
We will need this extended version later in the paper.
\end{remark}

\begin{proof}[Proof of Theorem~\ref{T1.5}] Let $z=0$ in \eqref{3.5} and use \eqref{1.23}.
\end{proof}

Killip--Nenciu \cite{KN2} obtain a limiting formula that involves sums of determinants,
but one can manipulate Heine's formula (see (1.5.80) of \cite{OPUC1}) to see they have
really found $\alpha_{\ell-k_{j-1}-1}(d\mu^{(K_j)})$. In fact, earlier in their proof they
essentially do an inverse of this process.

\section{Asymptotics of Real Schur Flows} \lb{s4}

In this section---the main one from the point of view of \eqref{1.1}---we will prove
Theorem~\ref{T1.1} and Proposition~\ref{P1.2}. The central object will be the nontrivial
probability measure $d\mu$ with
\begin{equation} \lb{4.1}
\alpha_n(d\mu) =\alpha_n (t=0)
\end{equation}
In terms of $d\mu$, we will be able to specify which case of Theorem~\ref{T1.1} holds.
Define
\begin{equation} \lb{4.2}
\Theta (d\mu) = \min \{\abs{\theta} \mid e^{i\theta}\in\sigma_\ess (d\mu)\}
\end{equation}
The central role of such extreme points of $\sigma_\ess$ for Toda flows was understood
by Deift--Li--Tomei \cite{DLT85}. Basically, if $x_N$ is interpreted as $\lim_{n\to\infty}
x_n$ when $N=\infty$, we will have
\begin{equation} \lb{4.3}
x_N=\cos(\Theta(d\mu))
\end{equation}

We begin by analyzing the case $\Theta =0$.

\begin{theorem}[Case (i) of Theorem~\ref{T1.1}]\lb{T4.1} Suppose $1\in\sigma_\ess
(d\mu)$. Then for all $n$ and $j=1, \dots, n$,
\begin{equation} \lb{4.4}
z_j^{(n)}(t)\to 1
\end{equation}
and
\begin{equation} \lb{4.5}
\alpha_n(t)\to (-1)^n
\end{equation}
\end{theorem}

\begin{proof} \eqref{4.5} follows from \eqref{4.4} and \eqref{1.23}. Let $z_1^{(n)}(t),
\dots, z_n^{(n)}(t)$ be the zeros of $\Phi_n (x;d\mu_t)$. By \eqref{3.4b},
\begin{equation} \lb{4.6}
z-z_j^{(n)}(t) = \Phi_1 \biggl( z; N_t^{-1} \prod_{k\neq j} \, \abs{z-z_k^{(n)}(t)}^2\, d\mu \biggr)
\end{equation}
where $N_t$ is a normalization. By Lemma~\ref{L4.2} below, the measure on the right converges to
$\delta_{z=0}$, so by continuity of $\alpha_n (\cdot)$ under weak convergence (with the critical
addendum that if the limit has $k$ pure points, it only holds for $n=0,1,\dots,k-1$),
\[
\alpha_0 \biggl( N_t^{-1} \prod_{k\neq j}\, \abs{z-z_k^{(n)}(t)}^2\, d\mu\biggr) \to 1
\]
so, since
\[
\Phi_1 =z-\bar\alpha_0
\]
we conclude
\begin{equation} \lb{4.7}
z_j^{(n)}(t) = \bar\alpha_0 \to 1
\end{equation}
Since $j$ is arbitrary, we have proven \eqref{4.4}.
\end{proof}

\begin{remark} For this case, where we only need information of $\Phi_1$, one can
use \eqref{3.2} instead of \eqref{3.4b}. By \eqref{3.2},
\[
z_j^{(n)} = \f{\int z\prod_{k\neq j}\abs{z-z_k^{(n)}}^2\, d\mu}
{\int \prod_{k\neq j}\abs{z-z_k^{(n)}}^2\, d\mu} \to 1
\]
by Lemma~\ref{L4.2}. We used \eqref{3.4b} since it is needed for the later arguments.
\end{remark}

\begin{lemma}\lb{L4.2} For any nontrivial probability measure $d\mu$ on $\partial\bbD$ with
$1\in\sigma_\ess(d\mu)$ and any $w_1(t), \dots, w_\ell (t)\in\ol{\bbD}$, we have
\begin{equation} \lb{4.9}
N_t^{-1} \prod_{j=1}^\ell\, \abs{z-w_j(t)}^2\, d\mu_t \to \delta_{z=1}
\end{equation}
as $t\to\infty$. Here $N_t =\int \prod_{j=1}^\ell \abs{z-w_j(t)}^2\, d\mu_t$.
\end{lemma}

\begin{proof} The idea is that in $\mu_t$, points near zero have much stronger weight than
fixed intervals away from zeros. The $\abs{z-w_\ell(t)}^2$ factors can overcome that
difference (and, as we have seen in the finite case, do if $d\mu$ has an isolated pure
point at $z=1$), but to do this, the $w_\ell(t)$ have to be exponentially close to
the points they mask. Thus, the finite number, $\ell$, of zeros can mask only an
exponentially small piece of the part of $d\mu_t$ near $z=1$, and since $1\in\sigma_\ess
(d\mu)$, there are always unmasked pieces.

To be explicit, let $d\ti\mu_t$ be the measure on the left side of \eqref{4.9}. If we prove
for each $\theta_0\in (0,\pi)$,
\begin{equation} \lb{4.09a}
\ti\mu_t(\{e^{i\eta}\mid \theta_0 <\eta<2\pi-\theta_0\})\to 0
\end{equation}
then, by compactness, \eqref{4.9} holds. Given $\theta_0$, since $1\in\sigma_\ess (d\mu)$, we
can find
\begin{equation} \lb{4.10}
\theta_0 >\theta_1 > \varphi_1 >\theta_2 > \varphi_2 >\cdots > \theta_{\ell+1} >
\varphi_{\ell+1} >0
\end{equation}
so that
\begin{equation} \lb{4.11}
\mu(\{e^{i\eta}\mid \theta_j >\eta >\varphi_j\})>0
\end{equation}
for each $j=1, \dots, \ell+1$. Call the set in \eqref{4.11}, $I_j$.

Let
\begin{equation} \lb{4.12}
Q_j(t)=\min_{\eta\in I_j}\, \prod_{m=1}^\ell \, \abs{e^{i\eta}-w_m(t)}^2
\end{equation}
Then, for each $j$,
\begin{equation} \lb{4.13}
\text{LHS of \eqref{4.09a}} \leq e^{2t(\cos\theta_0-\cos\theta_j)} \mu(I_j)^{-1}
Q_j(t)
\end{equation}
The right side of \eqref{4.13} goes to zero unless $Q_j(t)$ goes to zero as fast as
$e^{2t(\cos\theta_0 - \cos\theta_j)}$ and this can only happen if at least one $w_j(t)$
is within $e^{2t(\cos\theta_0-\cos\theta_j)/\ell}$ of $I_j$. Since
there are $\ell+1$ intervals a finite distance from each other and only $\ell$ zeros,
for all large $t$, at least one of the RHS of \eqref{4.13} (for $j=1, \dots, \ell+1$)
goes to zero, proving \eqref{4.09a}.
\end{proof}

The exact same argument proves:

\begin{lemma}\lb{L4.3} Let $d\mu$ be a nontrivial probability measure on $\partial\bbD$
invariant under $z\to\bar z$. Suppose $\Theta (d\mu)=\theta_\infty >0$ and $\{e^{i\eta}
\mid \abs{\eta}<\theta_\infty\}$ has finitely many pure points of $d\mu \colon
\{e^{i\theta_k}\}_{k=1}^K$ {\rm{(}}either one $\theta_k=0$ and $K$ is odd with $\pm
\theta_k$ terms or no $\theta_k$ is zero and $K$ is even with $\pm\theta_k$ terms{\rm{)}}.
Let $\ell\geq K$ and $w_1(t), \dots, w_\ell(t)$ be a conjugation-invariant set of
points in $\bbC$ so that for $j=1, \dots, K$,
\begin{equation} \lb{4.14}
\abs{e^{i\theta_j}-w_j(t)}^2 \leq Ce^{(2+\veps)t[\cos(\Theta)-\cos(\theta_j)]}
\end{equation}
for some $\veps >0$. Then
\begin{equation} \lb{4.15}
N_t^{-1} \prod_{j=1}^\ell\, \abs{z-w_j(t)}^2\, d\mu_t \to \tfrac12 (\delta_{z=e^{i\Theta}}
+\delta_{z=e^{-i\Theta}})
\end{equation}
\end{lemma}

\begin{remark} To get $\f12$ on the right in \eqref{4.15}, we use the fact that since
$d\mu_t$ and $\{w_j(t)\}$ are conjugation symmetric, the limit which lives on
$\{e^{\pm i\Theta}\}$ must also be conjugation symmetric.
\end{remark}

With this lemma, we can prove

\begin{theorem}\lb{T4.4} Let $d\mu$ be a nontrivial conjugation-symmetric probability
measure on $\partial\bbD$ with $\Theta(d\mu)=\theta_\infty >0$ and suppose there are only
finitely many points $\{e^{i\theta_k}\}_{k=1}^K$ in $\{e^{i\eta}\mid\abs{\eta}<\theta_\infty\}$
in the support of $d\mu$. Let $\alpha_n(t)$ solve \eqref{1.1} with
\begin{equation} \lb{4.16}
\alpha_n(0)=\alpha_n(d\mu)
\end{equation}

If $K=2m+1$ {\rm{(}}i.e., $\theta_1=0${\rm{)}} and $0<\theta_2 <\cdots < \theta_{m+1}$ and
$\theta_{m+2}=-\theta_{m+1}, \theta_{m+3}=-\theta_m, \dots, \theta_{2m +1} = -\theta_2$, then
\eqref{1.5} holds with $N=m+1$ and
\begin{align}
x_j&= \cos (\theta_{j+1}) \qquad j=1, \dots, m \lb{4.17} \\
x_N &= \cos (\theta_\infty) \lb{4.18}
\end{align}

If $K=2m$ and $0<\theta_1 < \cdots < \theta_m$, and $\theta_{m+1}=-\theta_m, \dots,
\theta_{2m}=-\theta_1$, then \eqref{1.6} holds with $N=m+1$ and
\begin{equation} \lb{4.19}
x_j = \cos(\theta_j) \qquad j=1, \dots, m
\end{equation}
and $x_N$ given by \eqref{4.18}.
\end{theorem}

\begin{proof} We will prove the case $K=2m+1$. The other case is essentially identical.
We essentially have a one-element $K_1$-group $\{z=1\}$ and $m$ two-element $K$-groups,
$K_2, \dots, K_{m+1}$, with $K_j = \{e^{\pm i\theta_j}\}$. The analysis of
Section~\ref{s3} (see the remark following Theorem~\ref{T3.4}) works for
$\{\alpha_n(t)\}_{n=0}^{2m}$ and proves that for any $\ell >2m+1$,
$\Phi_\ell (z;d\mu)$ has zeros exponentially close to $\{e^{i\theta_k}\}_{k=1}^K$
in the sense of \eqref{4.14} (indeed, one can take $2+\veps =4$).

Thus for any $\ell\geq 0$,
\begin{equation}\lb{4.21a}
N_t^{-1} \prod_{j=1}^{2m+1} \, \abs{z-z_j^{(2m+1+\ell)}}^2\, d\mu_t \to
\tfrac12\, (\delta_{z=e^{i\Theta}} + \delta_{z=e^{-i\Theta}})\equiv d\eta
\end{equation}
so by \eqref{3.4b} for $\ell=1,2$,
\begin{equation}\lb{4.21b}
\Phi_{2m+1+\ell}(z)\to \prod_{j=1}^{2m+1} (z-z_j) \Phi_\ell (z;d\eta)
\end{equation}

The two-point measure has $\alpha_0=\cos(\Theta)$ and $\alpha_1 =-1$, which proves
the formula for $\alpha_{2m+1}(d\mu_t)$ and $\alpha_{2m+2}(d\mu_t)$.

By Lemma~\ref{L4.4A} below and the argument in the first paragraph, we know that for
any $\ell\geq 0$, that $2m+3$ zeros of $\Phi_{2m+3+\ell}(t)$ approach $\{e^{i\theta_k}
\}_{k=1}^K \cup \{e^{\pm i\Theta}\}$. Repeating the argument above, we get $\alpha_{2m+3}
(d\mu_t)$ and $\alpha_{2m+4}(d\mu_t)$. Iterating, we get $\alpha_{2m+\ell}(d\mu)$
for all $\ell$.
\end{proof}

\begin{lemma} \lb{L4.4A} If $d\mu_t$ is a family of measures indexed by $t\in (0,\infty)$
and for some fixed $N$\!, there are $\{z_j^{(\infty)}\}_{j=1}^N\in\partial\bbD$ so that
the zeros $\{z_j^{(n)}(t)\}_{j=1}^N$ of $\Phi_N (z;d\mu_t)$ approach $\{z_j^{(\infty)}
\}_{j=1}^N$, then for any $\ell >N$\!, there are $N$ zeros of $\Phi_\ell (z;d\mu_t)$
which approach $\{z_j^{(\infty)}\}_{j=1}^N$.
\end{lemma}

\begin{proof} Since the coefficients of $\Phi_j (z;d\nu)$ are uniformly bounded by $2^j$
(uniformly in $d\nu$ by Szeg\H{o} recursion), $\Phi_\ell (z;d\mu)$ are uniformly bounded
analytic functions. So it suffices to show for each $j$, $\Phi_\ell (z_j^{(\infty)};d\mu_t)\to 0$.

Since $\Phi_N (z;d\mu_t)\to \prod_{j=1}^N (z-z_j^{(\infty)})$, we have
\begin{equation}\lb{4.21c}
\Phi_N^*(z;d\mu_t) \to \prod_{j=1}^N (1-z\bar z_j^{(\infty)}) = \prod_{j=1}^N
(-\bar z_j^{(\infty)}) \Phi_N (z;d\mu_t)
\end{equation}
Thus
\begin{equation}\lb{4.21d}
\Phi_N^*(z_j^{(\infty)}) \to 0 \qquad
\Phi_N (z_j^{(\infty)}) \to 0
\end{equation}
By Szeg\H{o} recursion,
\begin{equation}\lb{4.21e}
\Phi_{N+1}^*(z_j^{(\infty)}) \to 0 \qquad
\Phi_{N+1} (z_j^{(\infty)}) \to 0
\end{equation}
so by induction,
\begin{equation}\lb{4.21f}
\Phi_{N+m}^*(z_j^{(\infty)}) \to 0 \qquad
\Phi_{N+m} (z_j^{(\infty)}) \to 0
\end{equation}
for all $m$.
\end{proof}

We summarize in a strong version of Theorem~\ref{T1.1}:

\begin{theorem}\lb{T4.5} Suppose each $\alpha_n(0)\in (-1,1)$ and let $d\mu$ be the measure
with $\alpha(d\mu)=\alpha_n(0)$ {\rm{(}}which is conjugation-symmetric{\rm{)}}. Then
\begin{SL}
\item[{\rm{(i)}}] If $1\in\sigma_\ess(d\mu)$, \eqref{1.3} holds.

\item[{\rm{(ii)}}] If $\Theta(d\mu)>0$, $1\notin\supp(d\mu)$, and $\{e^{i\eta}\mid
\abs{\eta} < \Theta\}$ has $2m$ points, then \eqref{1.6} holds with $N=m+1$ and
$x_j$ is given by \eqref{4.19} and \eqref{4.18}.

\item[{\rm{(iii)}}] If $\Theta(d\mu)>0$, $1\notin\supp(d\mu)$, and $\{e^{i\eta}\mid
\abs{\eta}<\Theta\}$ has an infinity of points, then \eqref{1.6} holds with $N=\infty$
and $x_j$ is given by \eqref{4.18} and $x_j\to\cos(\Theta)$ as $N\to\infty$.

\item[{\rm{(iv)}}] If $\Theta(d\mu) >0$, $1\in\supp(d\mu)$, and $\{e^{i\eta}\mid
\abs{\eta}<\Theta\}$ has $2m+1$ points, then \eqref{1.5} holds, $N=m+1$, and
$\theta_j$ is given by \eqref{4.17} and \eqref{4.18}.

\item[{\rm{(v)}}] If $\Theta(d\mu) >0$, $1\in\supp(d\mu)$, and $\{e^{i\eta}\mid
\abs{\eta}<\Theta\}$ has an infinity of points, then \eqref{1.5} holds with
$N=\infty$ and $x_j$ is given by \eqref{4.17} and $x_j\to\cos(\Theta)$ as $N\to\infty$.
\end{SL}
\end{theorem}

\begin{proof} (i)--(iii) are proven in Theorems~\ref{T4.1} and \ref{T4.4}. (iv)--(v) follow
from the method of Section~\ref{s3} with no need for analysis of the edge of the
essential spectrum.
\end{proof}

\begin{proof}[Proof of Proposition~\ref{P1.2}] (i)\ By Theorem~4.3.17 of \cite{OPUC1},
$\alpha_n(0)\to 0$ implies $\supp(d\mu)=\partial\bbD$ implies $1\in\sigma_\ess(d\mu)$.

\smallskip
(ii)\ By Theorem~4.2.11 of \cite{OPUC1}, if $\alpha_{n+1}\bar\alpha_n\to
-1$, $\sigma_\ess (d\mu)=\{1\}$, so $1\in\sigma_\ess (d\mu)$.

\smallskip
(iii),(iv)\ By Example~4.3.10 of \cite{OPUC1}, $\sigma_\ess (d\mu)=[\theta_0,
2\pi-\theta_0]$ where $\cos\theta_0=1-2a^2$. So the only issue is whether $1\in\supp
(d\mu)$ or not. Since $\alpha_n(0)\sim \prod_{j=0}^n (1-\alpha_j)$ (for $\alpha_j$ real),
we see that $\alpha_j\to a>0$ means $1\in\supp(d\mu)$ and $\alpha_j\to a<0$ means
$1\notin\supp(d\mu)$.

\smallskip
(v),(vi)\ By Theorem~4.2.11 of \cite{OPUC1}, $\sigma_\ess (d\mu) =\{-1\}$
if and only if $\alpha_{n+1}\bar\alpha_n\to 1$, which for $\alpha_n$ real means
$\alpha_n\to 1$ or $\alpha_n \to -1$. If $\alpha_n\to 1$ (resp.\ $-1$), $1\in\supp(d\mu)$
(resp.\ $1\notin\supp(d\mu)$) by the argument used for (iii) and (iv).
\end{proof}

\section{Higher-Order Asymptotics in Case (i)} \lb{s5}

For the special case $\alpha_n (t=0)\equiv 0$, Golinskii proved that
\begin{equation}\lb{5.1}
(-1)^n \alpha_n(t) =1 - \f{n+1}{4t} + o\biggl(\f{1}{t}\biggr)
\end{equation}
(by our method below or his method, one can see $o(1/t)$ is $O(1/t^2)$).
He has two proofs: one uses the difference equation (in $n$) obeyed by $\alpha_n(t)$
in this special case, and the other, some explicit formulae in terms of Bessel functions.
Here is another proof which does not depend on special features of $\alpha_n(0)\equiv 0$
but only depends on the form $d\mu$ near $\theta=0$:

\begin{proposition} If $\alpha_n(0)$ is real and
\begin{equation}\lb{5.2}
\sum_{n=0}^\infty \, \abs{\alpha_n (t=0)} <\infty
\end{equation}
then \eqref{5.1} holds.
\end{proposition}

\begin{remark} All the proof requires is that $d\mu =w(\theta)\f{d\theta}{2\pi}
+d\mu_\s$ where $0\notin\supp(d\mu_\s)$ and $\theta=1$ is a Lebesgue point of $w$
with positive density. \eqref{5.2} is only used to prove that.
\end{remark}

\begin{proof} We use the following formula (see (1.5.80) and (1.5.88) of \cite{OPUC1}):
\begin{equation}\lb{5.2a}
(-1)^n \alpha_n (d\mu) = N^{-1} \int e^{-i(\theta_0 + \cdots + \theta_n)}
\prod_{0\leq j < k\leq n} \abs{e^{i\theta_j} - e^{i\theta_k}}^2
\prod_{j=0}^n d\mu(\theta_j)
\end{equation}
where
\[
N=\int \prod_{0\leq j < k\leq n} \abs{e^{i\theta_j} -e^{i\theta_k}}^2
\prod_{j=0}^n  d\mu(\theta_j)
\]

By Baxter's theorem, (see Theorem~5.2.1 of \cite{OPUC1}), \eqref{5.2} implies
\begin{equation}\lb{5.3}
d\mu_{t=0}(\theta) = f(\theta)\, \f{d\theta}{2\pi}
\end{equation}
where $f$ is continuous and nonvanishing.

By \eqref{5.2a},
\begin{equation}\lb{5.4}
\begin{split}
1 &- (-1)^n  \alpha_n (d\mu_t) \\
& = N_t^{-1} \int (-e^{-i(\theta_0 +\cdots +\theta_n)} +1)
\prod_{0\leq j<k<n}\, \abs{e^{i\theta_j} -e^{i\theta_k}}^2 \prod_{j=0}^n
e^{2t\cos\theta_j} f(\theta_j)\, \f{d\theta_j}{2\pi}
\end{split}
\end{equation}
For each fixed $\veps$, we can break the integral up into the region $\abs{\theta_j}<\veps$
for all $\veps$, and its complement. The integral over the complement is bounded in
absolute value by $Ce^{2t\cos(\veps)}$.

Consider next the integral over the region $\abs{\theta_j}<\veps$. By the $\theta_j
\to -\theta_j$ symmetry, we can replace $1-e^{-i(\theta_0+\cdots +\theta_n)}$ by
$1-\cos  (\theta_0 +\cdots +\theta_n)$ and so get a positive integrand. Picking
the contribution of the region
\begin{equation}\lb{5.5}
\f{2j\veps}{(n+1)^2} < \theta_j < \f{(2j+2)\veps}{(n+1)^2} \equiv \eta_j(\veps)
\end{equation}
we get a lower bound of the form
\[
C\veps^2 (\veps^2)^{n(n-1)/2} e^{2t(\sum_{j=0}^n \cos\eta_j)}(\veps)^n
\]
In this way, we see that for $\veps$ small, the ratio of the complement to the
remainder is
\[
O(\veps^{-p} e^{-tD\veps^2})
\]
for some $p,D>0$ (which are $n$ dependent). This goes to zero as $t\to\infty$ for
any fixed $n,\veps$, so in $N_t$ and the integral, we can restrict integrals
to $\abs{\theta_j}<\veps$ and make an arbitrarily small fractional error.

Once $\abs{\theta_j}<\veps$, we can replace $f(\theta)$ by $f(0)$, $e^{2t\cos\theta_j}$
by $e^{2t-t\theta_j^2}$, and $e^{i\theta_j}-e^{i\theta_k}$ by $\theta_j-\theta_k$
with fractional errors going to zero. We conclude that
\begin{equation}\lb{old5.4}
\begin{split}
1-&(-1)\alpha_n(d\mu_t) = N_t^{-1}\int \tfrac12\, (\theta_0 +\cdots + \theta_n)^2 \\
&\prod_{0\leq j < k\leq n} \abs{\theta_j-\theta_k}^2
e^{-t(\theta_0^2+\cdots +\theta_n^2)} \, \f{d\theta_0}{2\pi} \cdots \f{d\theta_n}{2\pi}
\, (1+o(1))
\end{split}
\end{equation}
where $N_t$ is the integral without the $\f12 (\theta_0 +\cdots + \theta_n)^2$.

To see this, we use positivity. Thus $f(0) >0$ and $f$ continuous means $\delta(\veps)
=\sup_{\abs{\theta} <\veps} \abs{\f{f(\theta)}{f(0)} -1} \to 0$ as $\veps\downarrow 0$,
so using
\[
f(0)(1-\delta) \leq f(\theta) \leq f(0)(1+\delta)
\]
in all places, we establish the claim about fractional errors in replacing $f(\theta)$ by
$f(0)$. Similarly,
\[
\lim_{\veps\downarrow 0}\, \sup_{\substack{ \abs{\theta_j} <\veps \\
\abs{\theta_k} <\veps}}\, \biggl| \f{e^{i\theta_j}-e^{i\theta_k}}{\theta_j-\theta_k} -1\biggr|
=0
\]
allowing the other fractional replacement. Because of $t$-dependence, the replacement is
more subtle in $e^{2t\cos(\theta_j)}$. We show there is $\delta(\veps)$ so $\abs{\delta}<\veps$
means
\[
e^{-(1+\delta (\veps))\theta^2 t} \leq e^{2t(1-\cos\theta)} \leq e^{-(1-\delta(\veps))\theta^2 t}
\]
and that allows, after calculations below, to show, as $\delta\downarrow 0$, the replacement is allowed.

In the integrals in \eqref{old5.4}, we can now take $\theta_j$ to run from $-\infty$
to $\infty$ for, by the same arguments as above, those integrals are dominated by
the region $\abs{\theta_j}<\veps$.

Now change variables from $\theta_0, \dots, \theta_n$ to
\begin{align*}
x_0 &= \theta_0 +\cdots + \theta_{n+1} \\
x_j &= \theta_j -\theta_0 \qquad j=1,\dots, n
\end{align*}
Since $e_0 =\delta_0 +\cdots +\delta_n$ is a vector of Euclidean length $\sqrt{n+1}$
and $\delta_j-\delta_0\perp e_0$, we see
\[
\theta_0^2+\cdots + \theta_n^2 =\f{x_0^2}{n+1}+ Q(x_1, \dots, x_n)
\]
where $Q$ is a positive quadratic form.

The integrands in \eqref{old5.4} factor into functions of $x_0$ or $(x_1, \dots, x_n)$
with identical integrands in $(x_1, \dots, x_n)$. Thus, those factors cancel and we are
left with $x_0$ integrals only, and we get:
\begin{align*}
\text{Integral on RHS of \eqref{old5.4}}
&= \f{\f12 \int y^2 e^{-ty^2/n+1}\, dy}{N} \\
&= \f{(\f{n+1}{2}) \int w^2 e^{-tw^2} \, dw}{N} \\
& = -\f{n+1}{2}\, \f{d}{dt} \biggl[ \log \int e^{-tw^2}\, dw\biggr] \\
&= \f{n+1}{4}\, \f{1}{t}
\end{align*}
proving \eqref{5.1}.
\end{proof}

Once we drop \eqref{5.2}, the higher-order asymptotics are not universal. For example, if
$d\mu = \f{d\theta}{4\pi} + \f12 \delta_{\theta=0}$ (for which $\alpha_n(0)= \f{1}{n+2}$;
see Example~1.6.3 of \cite{OPUC1}), then $1-\alpha_0(t)=ct^{-3/2}$. One can similarly get lots of
variation in asymptotics of $\alpha_n(t)$ staying within $\sum_{n=0}^\infty \abs{\alpha_n(0)}^2
<\infty$. In fact, one can arrange that $\log (1-\alpha_\bddot(t))/\log t$ does not have a
limit as $t\to\infty$.

\section{Pathologies for Complex Initial Conditions} \lb{s6}

Our goal in this section is to provide the example that proves Proposition~\ref{P1.3}. As
a warmup, we consider a measure on $[0,\infty)$:

\begin{example}\lb{E6.1} Let
\begin{equation}\lb{6.1}
x_k =\f{1}{k!}
\end{equation}
and
\begin{equation}\lb{6.2}
\rho_k = e^{-k^{3/2}}
\end{equation}
and
\begin{equation}\lb{6.3}
d\rho_t = \f{\sum_{k=1}^\infty e^{-tx_k} \rho_k \delta_{x_k}}
{\sum_{k=1}^\infty e^{-tx_k} \rho_k}
\end{equation}
Finally, let
\begin{equation}\lb{6.4}
t_n =n!
\end{equation}
Then, we claim that
\begin{equation}\lb{6.5}
\| d\rho_{t_n} -\delta_{x_n}\| = O(e^{-n^{1/2}})
\end{equation}
so $\rho_{t=n!}$ is concentrated more and more at single but variable points.

It is not hard to see that \eqref{6.5} is equivalent to
\begin{equation}\lb{6.6}
\f{\sum_{m\neq n} e^{-t_n x_m} \rho_m}{e^{-t_n x_n} \rho_n} =
O(e^{-n^{1/2}})
\end{equation}
or, since $t_n x_n=1$, to
\begin{equation}\lb{6.7}
\sum_{m\neq n} e^{-t_n x_m} \rho_m \rho_n^{-1} = O(e^{-n^{1/2}})
\end{equation}

Note first that
\begin{align*}
e^{-t_n x_{n-j}} &= e^{-n!/(n-j)!} \\
&=\begin{cases} e^{-n} & j=1 \\
O(e^{-n(n-1)}) & j\geq 2
\end{cases}
\end{align*}
where the $O(\dots)$ term is uniform in $j$ and $n$. On the other hand,
\[
\rho_{n-1} \rho_n^{-1} = e^{n^{3/2}-(n-1)^{3/2}} = e^{O(n^{1/2})}
\]
while
\[
\rho_{n-j} \rho_n^{-1} \leq \rho_n^{-1} = e^{n^{3/2}}
\]
for $j\geq 2$. It follows that
\[
\sum_{m<n} e^{-t_n x_m} \rho_m \rho_n^{-1} = e^{-n + O(n^{1/2})} +
ne^{-n^2 + O(n^{3/2})}
\]
is certainly $O(e^{-n^{1/2}})$.

For $m>n$, we only need $e^{-t_n x_{n+j}}\leq 1$ since
\begin{align*}
\rho_{n+j} \rho_n^{-1} &= e^{-(n+j)^{3/2}} e^{n^{3/2}} =
\exp \biggl( -\int_n^{n+j} \tfrac32\, x^{1/2}\, dx \biggr) \\
&= O(e^{-jn^{1/2}})
\end{align*}
so
\[
\sum_{m>n} e^{-t_n x_m} \rho_m \rho_n^{-1} = O\biggl(\,\sum_{j=1}^\infty
e^{-jn^{1/2}}\biggr) = O(e^{-n^{1/2}})
\]
This proves \eqref{6.7}, and so \eqref{6.5}. \qed
\end{example}

\begin{example}\lb{E6.2} We will re-use Example~\ref{E6.1} but with the wrinkle that
for odd $n$, we will put the points near $i$, and for even $n$, near $-i$.
Explicitly, define $z_n$ by
\begin{gather}
2\Real z_n = -\f{1}{n!} \lb{6.8} \\
\abs{z_n}=1 \lb{6.9} \\
(-1)^{n+1} \Ima z_n > 0 \lb{6.10}
\end{gather}
With
\[
\mu_n = \f{e^{-n^{3/2}}}{\sum_{j=1}^\infty e^{-j^{3/2}}}
\]
define
\begin{equation}\lb{6.11}
d\mu = \sum_{n=1}^\infty \mu_n \delta_{z_n}
\end{equation}
and $d\mu_t$ by \eqref{1.11}.

As in Example~\ref{E6.1},
\begin{equation}\lb{6.12}
d\mu_{t=n!} -\delta_{z_n} = O(e^{-n^{1/2}})
\end{equation}
so since
\begin{equation}\lb{6.13}
\ol{\alpha_0(t)} =\int z\, d\mu_t
\end{equation}
we have
\begin{equation}\lb{6.14}
\alpha_0(n!) -\bar z_n\to 0
\end{equation}
But $z_{2n}\to -i$, while $z_{2n+1}\to i$, so
\[
\lim_{n\to\infty}\, \alpha_0((2n)!) = i \qquad
\lim_{n\to\infty}\, \alpha_0 ((2n+1)!) = -i
\]
and $\lim\alpha_0(t)$ does not exist. This proves Proposition~\ref{P1.3}.
\qed
\end{example}

Since after $z_n$, the dominant weight at time $t_n$ is at $z_{n+1}$, it should
be possible to show, using the methods of Sections~\ref{s3} and \ref{s4}, that
$\alpha_1 (n!)- \bar z_{n+1}\to 0$ and $\alpha_j (n!)- \bar z_{n+j}\to 0$ so that no
$\alpha_j(t)$ has a limit.

\section{Schur Flows with Complex but Regular Initial Conditions} \lb{s7}

If one allows complex initial conditions, there are two cases where the pathology
of the previous section does not occur. As earlier, $\Theta(d\mu)$ is defined by
\eqref{4.2}.

\begin{theorem}\lb{T7.1} If
\begin{equation}\lb{7.1a}
S=\{\theta\mid\abs{\theta}<\Theta,\, e^{i\theta}\in\supp(d\mu)\}
\end{equation}
is infinite, the asymptotic dynamics is determined as follows. Number the points
in $\{\theta_j\}_{j=0}^\infty$ in $S$ so
\begin{equation}\lb{7.1b}
\abs{\theta_0} \leq \abs{\theta_1}\leq \abs{\theta_2}\leq \dots
\end{equation}
We have that
\begin{SL}
\item[{\rm{(i)}}] If $\abs{\theta_{j-1}} < \abs{\theta_j} < \abs{\theta_{j+1}}$, then
\begin{equation}\lb{7.1}
\alpha_j(t) \to (-1) \prod_{k=0}^j (-e^{-i\theta_k})
\end{equation}

\item[{\rm{(ii)}}] If $\abs{\theta_j}=\abs{\theta_{j+1}}$, then
\begin{equation}\lb{7.2}
\alpha_j(t) \to \biggl[\, \prod_{k=0}^{j-1} (-e^{-i\theta_k})\biggr]
[ae^{-i\theta_j}+(1-a)e^{i\theta_j}]
\end{equation}
where
\begin{equation}\lb{7.3}
a=\f{\beta_+}{\beta_+ +\beta_-}
\end{equation}
with
\begin{equation}\lb{7.4}
\beta_\pm = \prod_{k=0}^{j-1}\, \abs{e^{\pm i\theta_j} -e^{i\theta_k}}^2
\mu (\{e^{\pm i\theta_j}\})
\end{equation}
\end{SL}
\end{theorem}

\begin{theorem} \lb{T7.2} If the set $S$ of \eqref{7.1a} is finite and there is a
unique point, $e^{i\theta_\infty}$, in $\supp(\mu)$ with $\abs{\theta_\infty}
=\Theta$, then if the points $\{\theta_j\}_{j=0}^N$ in $S$ are labeled so \eqref{7.1b}
holds, then {\rm{(i), (ii)}} above hold for $j\leq N$, and for $j\geq N+1$,
\begin{equation}\lb{7.5}
\alpha_j(t)\to (-1) \biggl[\, \prod_{k=0}^N (-e^{-i\theta_k})\biggr]
(-e^{-i\theta_\infty})^{j-N}
\end{equation}
\end{theorem}

The proofs of these theorems are a simple modification of the arguments
in Sections~\ref{s3} and \ref{s4}. Essentially, so long as $\abs{\theta_j}
<\abs{\theta_{j+1}}$, the zeros of $\Phi_j (z;d\mu_j)$ approach
$\{e^{i\theta_k}\}_{k=0}^j$, and in the case of a nondegenerate bottom
of the essential spectrum, all the extra zeros go to $e^{i\theta_\infty}$.

We want to consider a case where the bottom of the essential spectrum
is degenerate, that is, $e^{\pm i\theta}\in\sigma_\ess (d\mu)$, but $\mu$
is regular near $\pm\Theta$ in a sense we will make precise. For simplicity,
we will suppose the set $S$ is empty and $e^{\pm i\Theta}$ are not pure
points. It is easy to handle $S$ finite or $\{e^{\pm i\Theta}\}$
eigenvalues, but it complicates the statements (and, of course,
$S$ infinite is already in Theorem~\ref{T7.1}).

\noindent{\bf Definition.} We say $\mu$ is {\it weakly regular\/} if there exists $a_\pm
\in (-1,\infty)$ so that
\begin{alignat}{2}
&\text{(i)} \quad && \mu (\{e^{i\theta}\mid\abs{\theta}\leq\Theta\}) =0 \lb{7.6} \\
&\text{(ii)} \quad && \lim_{\veps\downarrow 0}\, \f{\log(\mu \{e^{\pm i\theta}\mid
\Theta\leq \theta \leq \Theta +\veps\})}{\log \veps} =a_\pm \lb{7.7} \\
\intertext{If (ii) is replaced by}
&\text{(iii)} \quad && \lim_{\veps\downarrow 0}\, \mu (\{e^{i\theta}\mid \Theta \leq
\theta \leq \Theta +\veps\}) \veps^{-a_\pm} = C_\pm \in (0,\infty) \notag
\end{alignat}
we say $\mu$ is {\it strongly regular}.

\smallskip
Here are the theorems in this case:

\begin{theorem}\lb{T7.3} If $\mu$ is weakly regular with $a_+ - a_- \notin 2\bbZ$,
then the asymptotic dynamics is as follows. Suppose $n\geq 0$ is an integer so
$a_- \in (a_+ + 2n, a_+ + 2n+2)$ {\rm{(}}otherwise, if $a_- < a_+$, interchange them{\rm{)}}:
\begin{alignat*}{3}
&\text{\rm{(i)}} \quad && \alpha_j (t)\to -(-e^{-i\Theta})^j \qquad && j=0,1,\dots, n \\
&\text{\rm{(ii)}} \quad && \alpha_{n+2k}(t) \to -(-e^{-i\Theta})^n
\qquad && k=0,1, \dots \\
&\text{\rm{(iii)}} \quad && \alpha_{2n+2k+1}(t) \to -(-e^{-i\Theta})^{n-1} \qquad && k=0,1,\dots
\end{alignat*}
\end{theorem}

\begin{theorem}\lb{T7.4} If $\mu$ is strongly regular and $a_+ -a_-\in 2\bbZ$, say $a_-=
a_+ + 2n$ with $n\geq 0$, then
\begin{alignat*}{3}
&\text{\rm{(i)}} \quad && \alpha_j(t)\to -(-e^{-i\Theta})^j \qquad && j=0,\dots, n-1 \\
&\text{\rm{(ii)}} \quad && \alpha_{2n+2k-1}(t) \to -(-e^{-i\Theta})^{n-1} \qquad
&& k=0,1,\dots \\
&\text{\rm{(iii)}} \quad && \alpha_{n+2k}(t) \to (-e^{-i\Theta})^{n-1} (ae^{-i\Theta}
+(1-a)e^{i\Theta})
\end{alignat*}
with $a$ given by \eqref{7.3} with
\begin{equation}\lb{7.8}
\beta_+ = C_+ \qquad \beta_- = \abs{e^{i\Theta} -e^{-i\Theta}}^{2n-2} C_-
\end{equation}
\end{theorem}

The proofs here are simple modifications of the arguments in Section~\ref{s4}. In case
$a_+ -a_-\notin 2\bbZ$, $\mu_t$ approaches $\delta_{e^{i\Theta}}$, then $\abs{z-e^{i\Theta}}^2
\mu_t/N_t^{(1)}$ has $\delta_{e^{i\Theta}}$ as its limit if $\delta_{e^{i\Theta}}$ if
$a_- > a_+ +2$ or $\delta_{e^{-i\Theta}}$ if $a_- < a_+ +2$. We repeat this $n$ times, and
after that the limits alternate between $\delta_{e^{i\Theta}}$ and $\delta_{e^{-i\Theta}}$.
If $a_+ - a_- \in 2\bbZ$, we get measure modification by products of $\abs{z-z_j}^2$ which
are $a\delta_{e^{i\Theta}} + (1-a) \delta_{e^{-i\Theta}}$.


\section*{Appendix A: The Symes--Deift--Li--Tomei Representation of Dynamics} \lb{AppA}

\renewcommand{\theequation}{A.\arabic{equation}}
\renewcommand{\thetheorem}{A.\arabic{theorem}}
\setcounter{theorem}{0}
\setcounter{equation}{0}

In this appendix, we will find an operator formulation of the flow $\Sigma_t^G$ of \eqref{1.13}.
This formulation for the Toda flow was discovered by Symes \cite{Symes} and
then generalized to Jacobi analogs of $\Sigma_t^G$ by Deift--Li--Tomei \cite{DLT85}.
Killip--Nenciu \cite{KN2} discussed this for polynomial $G$ and finite CMV matrices.
We include this appendix for four reasons. First, we wish to show that one can
handle infinite CMV matrices. Second, while we regard this as a central result in
symplectic flows, it is somewhat hidden in \cite{KN2}. Third, the elementary identification
of the spectral measure below, while implicit in the earlier works, is not made explicit.
Finally, we want to note some aspects of the equivalence result (his Theorem~1) of
Golinskii \cite{Golppt} without extensive calculation.

The $QR$ algorithm is critical to this appendix. We will consider bounded operators on
$\calH\equiv\ell^2 (\{0,1,2,\dots\})$ which are therefore given by semi-infinite matrices.
We use $\{\delta_n\}_{n=0}^\infty$ for the canonical basis of $\calH$.

A bounded operator, $B$, is called positive upper triangular if and only if $B_{rs}=0$
for $r>s$ and $B_{rr}>0$. The set of all such operators will be denoted $\calR$.
The following is well known. We sketch the proof to emphasize the final formulae:

\begin{lemma}\lb{LA.1.1} Let $A$ be an invertible bounded operator on $\calH$. Then
\begin{equation} \lb{A.1.1}
A=QR
\end{equation}
with $Q$ unitary and $R\in\calR$. This decomposition is unique. Moreover,
\begin{equation} \lb{A.1.2}
Q\delta_0 = \f{A\delta_0}{\|A\delta_0\|}
\end{equation}
\end{lemma}

\begin{proof} Uniqueness is immediate since $Q_1 R_1=Q_2 R_2$ implies $Q_2^{-1} Q_1 =
R_2 R_1^{-1}$ and $B\in\calR$ and unitary implies $B=1$.

Let $e_0, e_1, \dots$ be the set obtained by applying Gram--Schmidt to $A\delta_1,
A\delta_2, \dots$. Note that
\begin{equation} \lb{A.1.3}
e_0 = \f{A\delta_0}{\|A\delta_0\|}
\end{equation}

Because $A$ is invertible, $\{e_n\}_{n=0}^\infty$ is a basis and
\begin{equation} \lb{A.1.4}
A\delta_j = \sum_{k=0}^n r_{kj} e_j
\end{equation}
and
\begin{equation} \lb{A.1.4a}
r_{kk} >0
\end{equation}
by the Gram--Schmidt construction.

Let $Q$ be defined by
\begin{equation} \lb{A.1.5}
Q\delta_j=e_j
\end{equation}
so \eqref{A.1.3} becomes \eqref{A.1.2}. By \eqref{A.1.4} and \eqref{A.1.4a},
\begin{equation} \lb{A.1.6}
Q^{-1} A\equiv R
\end{equation}
lies in $\calR$, and clearly, $A=QR$.
\end{proof}

We freely use the CMV matrix, $\calC$, and alternate CMV matrix, $\ti\calC$, discussed
in \cite{OPUC1} and \cite{CMVrvw}. Here is the main result of this appendix:

\begin{theorem}\lb{TA.1.2} Let $\calC$ be a CMV matrix associated to the measure
$d\mu$. Let $G$ be a real-valued function in $L^\infty (\partial\bbD, d\mu)$. Define
$Q_t, R_t$ by
\begin{equation} \lb{A.1.7}
\exp (\tfrac12\, t G(\calC)) = Q_t R_t
\end{equation}
using the $QR$ algorithm \eqref{A.1.1}. Define $\calC_t$ by
\begin{equation} \lb{A.1.8}
\calC_t = Q_t^{-1} \calC Q_t
\end{equation}
Then $\calC_t$ is the CMV matrix of the measure $d\mu_t$ given by \eqref{1.13}.
\end{theorem}

Our proof is related to that of Killip--Nenciu \cite{KN2}, which in turn is a
CMV analog of the results of Deift--Li--Tomei \cite{DLT85} for Jacobi matrices.
For notational simplicity, we will deal with nontrivial $d\mu$. \cite{KN2}
handles the case where $d\mu$ has finite support. We need:

\smallskip
\noindent{\bf Definition.} A matrix, $M$\!, on $\calH$ is said to have CMV shape if and
only if
\begin{SL}
\item[(i)] $M$ is five-diagonal, that is, $M_{jk}=0$ if $\abs{j-k} >2$
\item[(ii)] $M_{2n, 2n+2} >0$, $n=0,1,2, \dots$
\item[(iii)] $M_{2n+1, 2n+3} =0$, $n=0,1,2,\dots$
\item[(iv)] $M_{2n+3, 2n+1} >0$, $n=0,1,2,\dots$
\item[(v)] $M_{2n+2, 2n} =0$, $n=0,1,2,\dots$
\item[(vi)] $M_{10} >0$
\end{SL}
We say $M$ has alternate CMV shape if $M^t$ has CMV shape.

\smallskip
For finite matrices, the following is a result of \cite{KN2}:

\begin{proposition}\lb{PA.1.3} A unitary matrix, $M$\!, has CMV shape if and only if
for some sequence of Verblunsky coefficients $\{\alpha_n\}_{n=0}^\infty \subset
\bbD^\infty$, $M=\calC(\{\alpha_n\}_{n=0}^\infty)$. It has alternate CMV shape
if and only if $M=\ti\calC (\{\alpha_n\}_{n=0}^\infty)$.
\end{proposition}

\begin{remark} Our proof differs from \cite{KN2} in that they use a Householder
algorithm and we use the simpler AMR factorization (see \cite{CMVrvw}).
\end{remark}

\begin{proof} That a CMV matrix has CMV shape follows from the form of $\calC$;
see (4.2.14) of \cite{OPUC1}.

For the converse, define $\alpha_0\in \bbD$ and $\rho_0 \in (0,1)$ by
$M\delta_0 = \binom{\bar\alpha_0}{\rho_0}$. By (v), $\rho_0 >0$ and, by
unitarity, $\abs{\alpha_0}^2 + \rho_0^2=1$. Let
\begin{equation} \lb{A.1.9}
\Theta(\alpha) = \begin{pmatrix}
\bar\alpha & \rho \\
\rho & -\alpha
\end{pmatrix}
\end{equation}
and consider
\begin{equation} \lb{A.1.10}
(\Theta (\alpha_0)\oplus\bdone)^{-1} M
\end{equation}
It is clearly unitary and has $1$ in the $11$ corner, and by $\rho_0 >0$ and the
definition of CMV shape, it is of the form $\bdone_{1\times 1}\oplus M_1$ where
$M_1$ is of alternate CMV shape. Thus
\begin{equation} \lb{A.1.11}
M=(\Theta(\alpha_0)\oplus\bdone)(\bdone_{1\times 1} \oplus M_1)
\end{equation}

Applying this to get $M_1^t$, we see
\begin{equation} \lb{A.1.12}
M=(\Theta (\alpha_0)\oplus\bdone)(\bdone_{2\times 2}\oplus M_2)
(\bdone_{1\times 1}\oplus \Theta (\alpha_1) \oplus\bdone)
\end{equation}
where $M_2$ is of CMV shape. Iterating this $n$ times,
\begin{equation} \lb{A.1.13}
\begin{split}
M &=(\Theta (\alpha_0) \oplus\cdots\oplus\Theta(\alpha_{2n-2})\oplus\bdone)
(\bdone_{2n\times 2n}\oplus M_{2n})\\
&\qquad\quad (\bdone_{1\times 1}\oplus\cdots\oplus\Theta(\alpha_1) \oplus \Theta(\alpha_2)
\oplus \cdots \oplus \Theta(\alpha_{2n-1})\oplus\bdone)
\end{split}
\end{equation}
where $M_{2n}$ is of CMV shape. Taking $n\to\infty$ and taking strong limits,
we see that $M=\calC (\{\alpha_n\}_{n=0}^\infty)$.
\end{proof}

\begin{proof}[Proof of Theorem~\ref{TA.1.2}] We can write
\begin{align}
\calC_t &= R_t (Q_t R_t)^{-1} \calC Q_t R_t R_t^{-1} \notag \\
&= R_t \calC R_t^{-1} \lb{A.1.13a}
\end{align}
since $e^{tG(\calC)/2}$ commutes with $\calC$. Since $R_t$ and $R_t^{-1}$ are
in $\calR$, it is easy to see that $(\calC_t)_{jk}=0$ if $k>j+2$ and the
conditions (ii) and (iii) of the definition of CMV shape hold.

On the other hand, since $e^{tG/2}$ is selfadjoint and $Q$ unitary:
\begin{equation} \lb{A.1.14}
e^{tG(\calC)/2} = (e^{tG(\calC)/2})^* = R_t^* Q_t^* = R_t^* Q_t^{-1}
\end{equation}
so
\begin{align}
\calC_t &= (R_t^*)^{-1} (R_t^* Q_t^{-1}) \calC (R_t^* Q_t^{-1})^{-1} R_t^* \notag \\
&= (R_t^*)^{-1} \calC R_t^* \lb{A.1.15}
\end{align}
so, since $R_t^*$ is lower triangular and positive on diagonal, $(\calC_t)_{jk}=0$
if $k<j-2$ and conditions (ii), (v), and (vi) hold.

Thus, $\calC_t$ has CMV shape, and so is a CMV matrix by Proposition~\ref{A.1.3}.

The spectral measure of $\calC_t$ and the vector $\delta_0$ is that of $\calC$ and
$Q_t \delta_0$ which, by \eqref{A.1.3}, is that of $\calC$ and $e^{tG(\calC)/2}\delta_0
/ \| e^{tG(\calC)/2}\delta_0\|$, which is $e^{tG(\theta)}
d\mu(\theta)/\int e^{tG(\theta)} d\mu(\theta)$.
\end{proof}

Finally, when $G$ is a Laurent polynomial, we want to discuss the associated
difference equation and the associated Lax form and, in particular, show that
\eqref{1.11} solves \eqref{1.1}. We will not be explicit about uniqueness, but it
is not hard to prove \eqref{1.1} has a unique solution. First, following the
Deift--Li--Tomei \cite{DLT85} calculation for the Toda analog:

\begin{proposition}\lb{PA.1.4} Define $\pi$ on selfadjoint matrices on $\calH$ by
\begin{equation} \lb{A.1.17}
\pi (A)_{jk} = \begin{cases}
A_{jk} & j<k \\
-A_{jk} & j>k \\
0 & j=k
\end{cases}
\end{equation}
so $\pi(A)$ is skew-adjoint. Then the $\calC_t$ of \eqref{A.1.8} is strongly $C^1$ and
obeys
\begin{equation} \lb{A.1.18}
\calC_t = [B_t, \calC_t]
\end{equation}
where
\begin{equation} \lb{A.1.19}
B_t = \pi (\tfrac12 G(\calC_t))
\end{equation}
\end{proposition}

\begin{example} For \eqref{1.11}, $G(\theta) =2\cos\theta$, $G(\calC_t)=\calC_t +
\calC_t^{-1}$ (under $\calC\leftrightarrow e^{i\theta}$) and
\begin{equation} \lb{A.1.20}
B_t = \tfrac12 \, [(\calC_t + \calC_t^{-1})_+ - (\calC_t + \calC_t^{-1})_- ]
\end{equation}
with $(\cdot)_+$ the part of the $(\cdot)$ above the diagonal and $(\cdot)_-$
below. \eqref{A.1.18}--\eqref{A.1.19} is, in this case, \eqref{1.21}--\eqref{1.22}
of \cite{Golppt}.
\qed
\end{example}

\begin{proof} Since Gram--Schmidt is an algebraic operation, on $\{\delta_n\}_{n=0}^\infty$,
$Q_t$ is strongly $C^1$, and so $R_t=Q_t^* \exp (\f12 tG(\calC))$ is strongly $C^1$.
Clearly, $Q_t R_t = e^{tG(\calC)/2}$ implies
\begin{equation} \lb{A.1.21}
\Dot{Q}_t R_t + Q_t \Dot{R}_t = \tfrac12 \, G(\calC) Q_t R_t
\end{equation}
or
\begin{equation} \lb{A.1.22}
Q_t^{-1} \Dot{Q}_t + \Dot{R}_t R_t^{-1} = \tfrac12\, Q_t^{-1} G(\calC) Q_t =
\tfrac12\, G(\calC_t)
\end{equation}

Since $R_t\in\calR$, $\Dot{R}_t R_t^{-1}$ is upper triangular and real on the diagonal.
Since $Q_t$ is unitary, $Q_t^{-1} \Dot{Q}_t$ is skew-Hermitian. It vanishes on
diagonal since both $\Dot{R}_t R_t^{-1}$ and $\f12 G(\calC_t)$ are real there and
skew-Hermitian matrices are pure imaginary on diagonals. Since $\Dot{R}_t R_t^{-1}$ is
upper triangular,
\begin{equation} \lb{A.1.23}
[Q_t^{-1} \Dot{Q}_t]_{jk} = [\tfrac12\, G(\calC_t)]_{jk}
\end{equation}
for $j>k$. Since $Q_t^{-1} \Dot{Q}_t$ is skew-Hermitian, vanishes on diagonal, and
$G(\calC_t)$ is Hermitian, we obtain
\begin{equation} \lb{A.1.24}
Q_t^{-1} \Dot{Q}_t = -B_t
\end{equation}
where $B_t$ is given by \eqref{A.1.19}.

Differentiating \eqref{A.1.8},
\begin{equation} \lb{A.1.25}
\Dot{\calC_t} = Q_t^{-1} [-\Dot{Q}_t] Q_t^{-1} \calC Q_t + Q_t^{-1} \calC \Dot{Q}_t
\end{equation}
Inserting $Q_t Q_t^{-1}$ before the final $\Dot{Q}_t$, we see
\[
\Dot{\calC}_t = B_t \calC_t - \calC_t B_t
\]
which is \eqref{A.1.18}.
\end{proof}

Finally, we note the following that can be obtained by taking limits of \cite{KN2} and
is discussed in \cite{GNprep,CSprep}. Let $H$ be a real Laurent polynomial, that is,
$\sum_{k=-n}^n c_k e^{ik\theta}$ where $c_{-k} =\bar c_k$ and let
\begin{equation} \lb{A.1.26a}
G(\theta) = \f{dH(\theta)}{d\theta}
\end{equation}
which is also a real Laurent polynomial. Define
\begin{align}
t_H (\{\alpha_j\}_{j=0}^\infty &= \text{``Tr''} (H(\calC(\alpha))) \lb{A.1.26} \\
&\equiv \sum_{k=-n}^n c_k \text{``Tr''} (\calC^k) \lb{A.1.27}
\end{align}
where ``Tr'' is a formal sum. While $t_H$ is a formal infinite sum, $\partial t_H/
\partial\bar\alpha_j$ is well defined since only finitely many terms depend on
$\bar\alpha_j$. Here is what is proven in \cite{KN2, GNprep,CSprep}:

\begin{proposition}\lb{PA.1.5} Let $H$ be a Laurent polynomial and $G$ given by
\eqref{A.1.26a}. Then $d\mu_t$ given by \eqref{1.11} solves
\begin{equation} \lb{A.1.28}
\Dot{\alpha}_j = i\rho_j^2 \, \f{\partial}{\partial\bar\alpha_j}\, (t_H(\alpha))
\end{equation}
\end{proposition}

\begin{example} If $H(e^{i\theta}) = 2\sin\theta$ so $G(e^{i\theta})=2\cos\theta$, then
\begin{align*}
t_H(\alpha) &= \text{``Tr''} (\tfrac12\, (\calC-\calC^*)) \\
&= i^{-1} \biggl( \, \sum_{j=0}^\infty\, [-\bar\alpha_j\alpha_{j-1} + \alpha_j
\bar\alpha_{j-1}]\biggr)
\end{align*}
(with $\alpha_{-1}\equiv -1$) and \eqref{A.1.28} becomes
\[
\Dot{\alpha}_j = \rho_j^2 (-\alpha_{j-1} + \alpha_{j+1})
\]
that is, \eqref{1.1}.
\qed
\end{example}

\section*{Appendix B: Zeros of OPUC near Isolated Points of the Spectrum} \lb{AppB}
\renewcommand{\theequation}{B.\arabic{equation}}
\renewcommand{\thetheorem}{B.\arabic{theorem}}
\setcounter{theorem}{0}
\setcounter{equation}{0}

In this appendix, we will prove a stronger result than Theorem~\ref{T3.1} using
a different proof from that of Denisov--Simon presented in \cite{OPUC1}. This will
use operator theory modeled on the following operator-theoretic proof of Fej\'er's theorem:

\begin{proposition}[well-known]\lb{PA.2.1} Let $A$ be a bounded normal operator. Then
for any unit vector, $\varphi$, $\langle \varphi, A\varphi\rangle$ lies in the
convex hull of the support of the spectrum of $A$.
\end{proposition}

\begin{proof} By the spectral theorem, $\varphi\in K$ a subspace (the cyclic subspace
generated by $A,A^*$ on $\varphi$), and there exists a probability measure $\mu$ on $\spec(A)$,
and $U\colon K\to L^2 (\bbC,d\mu)$ unitary, so $U\!AU^{-1}=$ multiplication by $z$,
$U\!A^*U^{-1} =$ multiplication by $\bar z$, and $U\varphi=1$. Thus
\[
\langle \varphi, A\varphi\rangle = \int z\,d\mu(z)
\]
clearly lies in the convex hull of $\spec(A)$.
\end{proof}

\begin{theorem}[Fej\'er's Theorem] \lb{TA.2.2} If $d\mu$ is a probability measure on $\bbC$
with
\begin{equation} \lb{A.2.1a}
\int \abs{z}^n\, d\mu <\infty
\end{equation}
for all $n$ and $\Phi_n(z)$ are the monic orthogonal polynomials for $d\mu$, then the zeros
of $\Phi_n$ lie in the convex hull of the support of $d\mu$.
\end{theorem}

\begin{proof} Let $P_n$ be the projection onto polynomials of degree up to $n-1$. Then
$P_n\eta =0$ for $\eta$ in $\ran P_{n+1}$ if and only if $\eta = c\Phi_n$. Thus, if
$Q\in \ran P$, then $P_n [(z-z_0)Q]$ if and only if $c\Phi_n(z)=(z-z_0)Q$ which happens
if and only if $\Phi_n (z_0)=0$. Thus, zeros of $\Phi_n$ are precisely eigenvalues of
$P_n M_z P_n \restriction \ran P_n$ where $M_z$ is multiplication by $z$.

If $\eta$ is the corresponding normalized eigenvector of $P_n M_z P_n$ in $\ran P_n$,
then $z_0=\langle \eta, P_n M_z P_n \eta\rangle = \langle \eta, z\eta\rangle$
lies in the convex hull of $\supp (d\mu)$ by Proposition~\ref{PA.2.1}.
\end{proof}

We are heading towards proving the following:

\begin{theorem}\lb{TA.2.3} Let $A$ be a normal operator and $z_0$ a simple eigenvalue
so that $z_0\notin\text{\rm{cvh}}(\spec(A)\backslash \{z_0\})\equiv C$. Here
\emph{cvh} is the ``convex hull of." Let $z_1\in C$ so that $\abs{z_0-z_1} =
\min\{\abs{w-z_0}\mid w\in C\}$. Let $P$ be an orthonormal projection and
$w_1, w_2$ two distinct eigenvalues of $P\!AP \restriction\ran P$. Then
\begin{equation} \lb{A.2.1}
\f{(w_1 - z_1)\, \bddot \, (z_0-z_1)}{\abs{z_1-z_0}^2} +
\f{(w_2 -z_1)\, \bddot\, (z_0 -z_1)}{\abs{z_1 - z_0}^2} \leq 1
\end{equation}
In particular,
\begin{equation} \lb{A.2.2}
\abs{w_1 - z_0} + \abs{w_2 - z_0} \geq \abs{z_1 - z_0}
\end{equation}
which implies there is at most one eigenvalue in
\begin{equation} \lb{A.2.3}
\{w\mid \abs{w-z_0} < \tfrac12\, \abs{z_1-z_0}\}
\end{equation}
\end{theorem}

\begin{remark} $\abs{z_1-z_0} = \dist (z_0, C)$.
\end{remark}

\begin{proof} Let $B=P\!AP$. Pick $\eta_1, \eta_2$ to be normalized eigenvalues of
$B$ and $B^*$, so
\begin{equation} \lb{A.2.4}
B\eta_1 = w_1 \eta_1 \qquad
B^* \eta_2 = \bar w_2 \eta_2
\end{equation}
and let $\varphi$ be the normalized simple eigenvector of $A$ with
\begin{equation} \lb{A.2.5}
A\varphi = z_0 \varphi \qquad
A^* \varphi = \bar z_0 \varphi
\end{equation}
(since $A$ is normal, if $\varphi$ obeys $A\varphi =z_0\varphi$, then $A^* \varphi$
also has $A(A^*\varphi) = z_0 A^*\varphi$ so $A^*\varphi =w\varphi$ and then
$\langle \varphi, A^*\varphi\rangle = \langle A\varphi,\varphi\rangle$ implies
$w=\bar z_0$).

By the spectral theorem for $A$,
\begin{equation} \lb{A.2.6}
w_1 =\langle \eta_1, A\eta_1\rangle=z_0 \abs{\langle \varphi,\eta_1\rangle}^2
+ (1-\abs{\langle\varphi,\eta_1\rangle}^2) x_1
\end{equation}
where $x_1\in C$. Here, if $Q$ is the projection onto multiples of $\varphi_1$,
then $x_1 = \langle(1-Q)\eta_1, A(1-Q\eta_1\rangle/ \|(1-Q)\eta_1 \|$.

By \eqref{A.2.6},
\begin{equation} \lb{A.2.7}
w_1 - z_1 = (1-\abs{\langle \varphi,\eta_1\rangle}^2)(x_1 -z_1) +
\abs{\langle \varphi, \eta_1\rangle}^2 (z_0 -z_1)
\end{equation}
Since $z_1$ minimizes $\abs{z_0-z_1}$ and $C$ is convex, $C$ is in the half-space orthogonal
to $z_0 - z_1$ not containing $z_0$ and so $(x_1 - z_1)\bddot (z_0-z_1) \leq 0$. Thus
\begin{equation} \lb{A.2.8}
(w_1-z_1)\, \bddot\, (z_0-z_1) \leq \abs{z_0-z_1}^2 \abs{\langle \varphi, \eta_1\rangle}^2\,
\end{equation}

Similarly, using
\[
w_2 = \ol{\langle\eta_2, A^*\eta_2\rangle} = \langle\eta_2, A\eta_2\rangle
\]
we see
\begin{equation} \lb{A.2.9}
(w_2-z_1) \bddot (z_0-z_1) \leq \abs{z_0-z_1}^2\, \abs{\langle \varphi,\eta_2\rangle}^2
\end{equation}

Thus, \eqref{A.2.1} is equivalent to
\begin{equation} \lb{A.2.10}
\abs{\langle \varphi,\eta_1\rangle}^2 + \abs{\langle\varphi, \eta_2\rangle}^2 \leq 1
\end{equation}

Next, note that
\begin{align*}
(w_1-w_2) \langle\eta_2,\eta_1\rangle &= \langle \eta_2, w_1 \eta_1\rangle -
\langle\bar w_2\eta_2, \eta_1\rangle \\
&= \langle B^*\eta_2, \eta_1\rangle - \langle \eta_2, B\eta_1\rangle =0
\end{align*}
so $\langle\eta_2, \eta_1\rangle=0$ and \eqref{A.2.6} is Bessel's inequality.
This proves \eqref{A.2.1}.

Obviously,
\begin{equation} \lb{A.2.11}
\f{(z_0-z_1)\, \bddot \, (z_0-z_1)}{\abs{z_0-z_1}^2} +
\f{(z_0-z_1)\, \bddot \, (z_0-z_1)}{\abs{z_0-z_1}^2} = 2
\end{equation}
so subtracting \eqref{A.2.1} from \eqref{A.2.7},
\[
\f{(z_0-w_1)\, \bddot \, (z_0-z_1)}{\abs{z_1 -z_0}^2} +
\f{(z_0-w_2)(z_0-z_1)}{\abs{z_1 -z_0}^2} \geq 1
\]
from which \eqref{A.2.2} follows by the Schwartz inequality.

That the set \eqref{A.2.3} can contain at most one $w$ is immediate from
\eqref{A.2.2} since $w_1, w_2\in$ the set \eqref{A.2.3} violates \eqref{A.2.2}.
\end{proof}

We can now improve the error of Denisov--Simon (who got $\delta/3$ where we get
$\delta/2$):

\begin{theorem} \lb{TA.2.4} Let $d\mu$ obey \eqref{A.2.1a} for all $n$ and let
$\Phi_n$ be the monic orthogonal polynomials. Suppose that $z_0$ is a pure point of
$d\mu$ and
\[
\delta = \dist (z_0, \text{\rm{cvh}} (\supp(d\mu)\backslash (\{z_0\}))) >0
\]
Then $\Phi_n$ has at most one zero in $\{z\mid \abs{z-z_0} < \delta/2\}$.
\end{theorem}

\begin{proof} Let $A$ be multiplication by $z$ on $L^2 (\bbC,d\mu)$ and $P_n$ be the
projection used in the proof of Theorem~\ref{TA.2.2}. Then Theorem~\ref{TA.2.3}
immediately implies the result.
\end{proof}

Specializing to OPUC:

\begin{theorem}\lb{TA.2.5} Let $d\mu$ be a probability measure on $\partial\bbD$
and $\{\Phi_n\}_{n=0}^N$ the monic OPUC {\rm{(}}if $d\mu$ is nontrivial, $N=\infty$;
if $d\mu$ has finite support, $N=\# (\supp(d\mu))${\rm{)}}. Let $z_0$ be an isolated
point of $\partial\bbD$ and $d=\dist (z_0, \supp(d\mu)\backslash \{z_0\})$. Then
for each fixed $n<N$, $\{z\mid\abs{z-z_0} < d^2/4\}$ has at most one zero of
$\Phi_n(z)$.
\end{theorem}

\begin{proof} As proven in \cite{OPUC1}, $\delta < d^2/2$.
\end{proof}

\section*{Acknowledgements}
It is a pleasure to thank Andrei Mart\'inez-Finkelshtein, Irina Nenciu, Paul Nevai,
and Vilmos Totik, and especially Leonid Golinskii for useful discussions and correspondence.

\medskip


\medskip

\end{document}